\documentclass[a4paper,DVI,10pt]{scrartcl}
\usepackage[utf8]{inputenc}
\usepackage{amssymb,amsmath,amsthm}
\usepackage{array, hhline}
\usepackage{tikz}
\usepackage{authblk}
\usepackage{subfigure}
\usepackage{array}
\usepackage{booktabs}
\usepackage{float}

\usetikzlibrary{plotmarks}
\renewcommand{\vec}{\boldsymbol}

\newcommand{\ud}{\,\mathrm{d}}
\newcommand{\R}{\mathbb{R}}
\renewcommand{\P}{\mathbb{P}}

\newtheorem{defi}{Definition}[section]
\newtheorem{theorem}[defi]{Theorem}

\newtheorem{rem}[defi]{Remark}
\newtheorem{cor}[defi]{Corollary}

\numberwithin{equation}{section}
\numberwithin{table}{section}
\numberwithin{figure}{section}

\newcommand{\Ls}{\ensuremath{\mathcal{L}}}
\newcommand{\Hs}{\ensuremath{\mathcal{H}}}
\newcommand{\Vs}{\ensuremath{\mathcal{V}}}

\newcommand{\Ps}{\ensuremath{\mathcal{P}}}

\newcommand{\T}{\ensuremath{\mathcal{T}}}

\newcommand{\J}{\ensuremath{\mathcal{J}}}
\newcommand{\St}{\ensuremath{\mathcal{S}}}
\newcommand{\DIV}{\nabla\cdot }
\newcommand{\eps}{\ensuremath{\varepsilon}}
\newcommand{\eref}[1]{(\ref{#1})}
\renewcommand{\vec}{\boldsymbol}
\newcommand\relph[1]{\mathrel{\phantom{#1}}}

\begin{document}

\title{\Large Goal-oriented a posteriori error control\\ for nonstationary 
convection-dominated\\ transport problems}
\author[K.\ Schwegler, M.\ P.\ Bruchh\"auser, M.\ Bause]
{\large K.\ Schwegler, M.\ P.\ Bruchh\"auser\thanks{bruchhaeuser@hsu-hh.de}, and 
M.\ Bause\thanks{bause@hsu-hh.de (corresponding author)}\\
{\small Helmut Schmidt University, Faculty of Mechanical Engineering,\\ 
Holstenhofweg 85, 220433 Hamburg, Germany}
}
\date{}
\maketitle

\begin{abstract}
{\bfseries Abstract.} 
The numerical approximation of convection-dominated problems continues to remain subject 
of strong interest. Families of stabilization techniques for finite element methods were 
developed in the past. Adaptive techniques based on a posteriori error estimates offer 
potential for further improvements. However, there is still a lack in robust a 
posteriori error estimates in natural norms of the discretizations. Here we combine the 
Dual Weighted Residual method for goal-oriented error control with stabilized finite 
element approximations. By a duality argument an error representation is derived on that 
a space-time adaptive approach is built. It differs from former works on the Dual 
Weighted Residual method. Numerical experiments illustrate that our schemes 
are capable to resolve layers and sharp fronts with high accuracy and to further reduce  
spurious oscillations of approximations. 
\end{abstract}


\section{Introduction}
\label{Sec:Intro}
In the last decades, since the pioneering works of the 1980's (cf., 
e.g., \cite{HMM86}), strong efforts were made in the development of accurate and 
efficient approximation schemes for convection-dominated flow and transport 
problems. For a review of fundamental concepts related to their analysis and 
approximation and a presentation of prominent robust numerical methods we refer 
to the monograph \cite{RST08}. Convection-dominated problems are of high 
practical interest. They arise in many branches of technology and, therefore, 
deservedly attracted substantial analysis. Applications can be found in fluid 
dynamics including turbulence modeling, electro-magnetism, semi-conductor 
devices, environmental and civil engineering as well as in chemical and 
biological scienes. If transport mechanisms are convection-dominated, solution 
profiles with sharp moving fronts, interior or boundary layers with complicated 
structures where important physical and chemical phenomena take place may occur. The 
development of numerical methods with the ability to capture the strong gradients of the 
exact solution without producing spurious oscillations or smearing effects continues to 
remain a challenging task. 

In the recent years a substantial progress has been made in the numerical approximation 
of convection-dominated problems, even though a real breakthrough is still missing. 
Numerous families of stabilization concepts were proposed and studied for various 
discretization techniques; cf.\ \cite{RST08}. Here we focus on finite element 
discretizations along with residual-based stabilizations. In particular, we use the 
streamline upwind Petrov--Galerkin (SUPG) method for our unsteady computations. For 
steady problems an additional shock-capturing stabilization is applied to further 
enhance the effect of the numerical method that is proposed in this work. Restricting 
ourselves to steady problems in some of our numerical experiments is done for the sake 
of simplicity. It is sufficient for illustrating the method's features we would like to 
investigate. In the literature shock-capturing stabilization is also often refered to as 
a spurious oscillations at layers diminishing (SOLD) method. The SUPG method reduces 
non-physical oscillations in streamline direction, whereas SOLD methods yield an 
additional stabilization in crosswind direction. For a review of prominent variants of 
SOLD methods and a competitive numerical investigation of the performance properties of 
SUPG and families of SOLD stabilizations we refer to, e.g., \cite{JS08}. Besides the 
class of these residual-based stabilization techniques, flux-corrected transport schemes 
are further addressed in \cite{JS08}. These techniques aim at a stabilization on the 
algebraic level; cf.\ \cite{KLT12}. In many works of the literature authors conclude that 
spurious oscillations in the numerical approximation of convection-dominated problems 
can be reduced by state-of-the-art stabilization techniques (cf., e.g., \cite{BS12}), but 
nevertheless the results are not satisfactory yet (cf.\ \cite{JS08}), in particular, if 
applications of practical interest and in three space dimensions are considered.        

Adaptive mesh generation based on an a posteriori error control is nowadays a well known 
technique to capture singular phenomena and sharp gradients of solutions to partial 
differential equations in numerical simulations. For a review of a posteriori error 
estimation techniques for finite element methods and automatic mesh generation we refer, 
for instance, to the monograph \cite{V13}. Even though interior and boundary layers that 
arise in applications of practical interest cannot be resolved completely by 
adaptive finite element meshes, at least in a reasonable computing time, a further 
improvement and gain in accuracy may nevertheless be expected by applying the concepts of 
automatic mesh generation to stabilized finite element approximations 
of convection-dominated problems. However, the design of an adaptive method requires the 
provision of an appropriate a posteriori error estimator. The derivation of such an error 
estimator for convection-dominated problems, that is robust with respect to the small 
perturbation parameter of the partial differential equation, is delicious 
and has borne out to be a considerable source of trouble. So far, the quality of 
adaptively refined grids is often not satisfactory yet. Further, only a few contributions 
have been published yet for convection-dominated problems and the considered type of 
discretizations. This observation even holds for stationary problems. For a deeper 
discussion and further references we refer to \cite{JN13} for the stationary case and to 
\cite{DEV13, FGJN15} for evolutionary problems. Usually, existing a posteriori 
error analyses for convection-dominated problems are either not robust with respect to 
the small perturbation parameter, embodied by increasing bounds for a vanishing 
perturbation parameter, or the a posteriori error estimates are not based on the natural 
norm of the discretization for that an a priori error analysis becomes feasible. For the 
stationary case, in \cite{V98} an a posteriori error bound is presented in the norm 
$(\|v\|_{\mathcal L^2(\Omega)}+ \varepsilon^{1/2}\|\nabla  v 
\|_{\mathcal L^2(\Omega)})^{1/2}$ with 
$\varepsilon$ being the small diffusion parameter. The bound is not robust in 
$\varepsilon$. On the other hand, in \cite{V05} the dual norm of the convective 
derivative is added to the energy norm to get a robust error estimate with respect to the 
small diffusion parameter. An extension of this approach to evolutionary problems is 
given in \cite{V05_2}. For evolutionary problems more recent robust a posteriori  
estimators measuring the error in a space-time mesh-dependent dual norm can be found in 
\cite{DEV13}. They are based on a space-time equilibrated flux reconstruction and are 
locally computable. The estimator by itself is local-in-time and local-in-space and does 
not depend on dual norms. However, dual error norms are usually beyond a reasonable 
physical interpretation and hard to compute such that they are of little interest in 
applications and difficult to use in studies of the experimental order of 
convergence. 

In \cite{JN13} a robust residual-based a posteriori error estimate in the SUPG norm is 
presented for stationary convection-diffusion equations. Its derivation uses 
variational multiscale theory. Upper and lower bounds are provided where the global upper 
bound relies on some hypotheses. A similar situation can be found in \cite{B09}. In 
\cite{JN13} it is argued and demonstrated by numerical experiments that the hypotheses 
are fulfilled and non restrictive in standard applications. The a priori error estimate in 
\cite{JN13} is based on different weights than other residual-based error estimators for 
convection-diffusion-reaction problems. However it is noted that the estimator performs 
well if the solution posses only one kind of singularity. Otherwise the non-robust 
residual-based estimator for the $\mathcal L^2$-norm should be prefered. In 
\cite{FGJN15} 
an adaptive SUPG method is proposed for time-dependent convection-diffusion problems where 
the SUPG solution is considered as a solution of a steady-state problem such that the 
error estimator of \cite{JN13} becomes applicable. However, the approach relies on the 
heuristic argument that a certain term is of higher order and thus becomes negligible. A 
validation of the assumption is given for one space dimension. In numerical calculations 
the robustness of the error estimator and its superiority over the adaptive approach that 
is presented in \cite{FGN11} and built upon heuristic error indicators is illustrated. A 
further non-robust a posteriori error of residual type estimator is presented in 
\cite{AV14}. Finally, we note that a posteriori error estimates are available for 
space-time finite element methods, cf.\ e.g. \cite{PP09}, and for Lagrange--Galerkin 
methods, cf.\ \cite{BC08,HLXZ14}, that are based on the method of characteristics and 
represent another class of prominent schemes for the approximation of evolutionary 
convection-dominated problems (cf.\ \cite{HS01}). A robust a priori error estimate for 
the Lagrange--Galerkin method with error constants depending only on norms of the data 
and not on (higher order) norms of the solution is presented in \cite{BK02}. The estimate 
is proved in a Lagrangian framework instead of using Eulerian coordinates as 
it is done in most of the error analyses for Lagrange-Galerkin methods. 

The Dual Weighted Residual method (or shortly DWR method) \cite{BR03} aims at the 
economical computation of arbitrary quantities of physical interest by properly adapting 
the computational mesh. Thus, the mesh adaption is based on the computation and control 
of a physically relevant goal quantity instead of the traditional energy-norm or the 
$\mathcal L^2$-norm. The DWR approach relies on a space-time variational formulation of 
the 
discrete problem and uses duality techniques to find a rigorous a posteriori error 
estimate. Such duality techniques are well known from a priori error analyses; cf., e.g., 
\cite{EJ91,EJ95,HR82}. The DWR approach has been applied to the numerical approximation 
of several classes of mathematical models based on partial differential equations, 
including fluid mechanics \cite{BR06}, wave propagation \cite{BGR10}, structural 
mechanics \cite{RS99}, fluid-structure interaction \cite{R12}, eigenvalue problems 
\cite{HR01}, optimization problems \cite{MB07} and, further, been applied to goal-oriented 
adaptive modeling \cite{BE03}. In the abstract DWR philosophy (cf.\ \cite{BR03} for 
details) an error representation for the considered goal quantity is derived by 
duality techniques at the beginning. This error identity cannot be evaluated directly, 
since it depends on the unknown solution of the ''linearized'' dual or adjoint problem 
that has to be solved numerically. We note that the dual problem is always a linear one, 
such that in the case of a nonlinear partial differential equation the numerical costs 
for solving the dual problem requires much less work. If the primal solution is obtained 
by a Newton iteration then solving the dual problem corresponds to one additional Newton 
iteration in each time step. From the error representation for the goal quantity 
localized error indicators can be derived (cf.\ \cite{BR03} for details), similarly to 
the traditional residual-based approach. However, the sharpness of the resulting a 
posteriori error estimate can not be guaranted anymore as soon as estimates are applied 
to the error identity.  

Even though the DWR approach has been applied to many classes of partial differential 
equations, our feeling is that its potential for the numerical approximation of 
convection-dominated problems and stabilized discretizations has not been completely 
understood and explored yet. In the application of the DWR method the efficient and fast 
solution of the auxiliary dual problem and the localization of the error estimator is an 
essential step in practice; cf.\, e.g., \cite{BBT11,BR03,BGR10,RW15,SV08}. The dual 
solution impacts the weights of the resulting error indicators. It is well known that the 
proper choice of the weights is crucial for the effectivity of the adaptation process. 
They should in particular measure the influence of a present cell on the requested goal 
quantity of interest. The approximation of the dual solution cannot be done in the finite 
element space of the primal problem since it would result in an useless vanishing 
approximation of the error quantity; cf.\ \cite{BR03}. Therefore, several techniques of 
approximation the dual solution efficiency were developed and proposed in the literature. 
Approximation by a higher-order method, approximation by a higher-order interpolation, 
approximation by difference quotients and approximation by local residual problems are 
known approaches \cite{BBT11,BGR10,BR03,SV08}. In particular, higher-order interpolation 
is applied often in the literature; cf., e.g., \cite{BGR10,SV08} .

In this work we combine the DWR approach with stabilized approximations of 
con\-vection-dominated problems. The adaptive mesh refinement process is directly based 
on the global error representation, up to higher order error contributions, without 
estimating the error terms further, i.e.\ without providing the usual upper bounds for 
the error in the goal quantity. Thereby we aim to reduce additional approximation errors 
and non-sharp estimates of standard error indicators and to apply the DWR approach as 
strictly as possible in order to get a reliable quantification and control of the error 
in the goal quantity instead of providing only error indicators for an economical mesh 
adaption strategy. In our approach the dual problem is solved by using higher oder 
finite element techniques. In numerical experiments we will illustrate the high impact of 
the proper choice of the weights and thereby of the dual solution on the mesh adaption 
process. Our overall motivation is to reduce sources of inaccuracies and non-sharp 
estimates in the a posteriori error control mechanism as far as possible in order to avoid 
numerical artefacs and a loss of quality in the approximation of the solution and the 
error control in regions with sharp fronts where highly sensitive solution profiles are 
present and interpolations are defective. We expect that this strategy allows us to 
exploit the full potential of the DWR method for the a posteriori quantification of 
discretization errors. This is in contrast to other works of the literature on the DWR 
method where much effort is put in the reduction of the computation costs for solving the 
dual problem. Thereby, the high impact of the dual solution on the error control and mesh 
generation process is not focused as strongly as in this work. Due to the specific 
character of convection-dominated problems we are convinced that the error control needs 
a particular care in regions with interior and boundary layers and in regions with sharp 
fronts in order to get an accurate quantification of numerical errors. For problems with 
simpler structures of solutions more economical approximations of the dual solution might 
be sufficient and appropriate.

This work is organized as follows. In Section \ref{Sec:pro}, we introduce our model 
problem together with some global assumptions and our notation. Further we present the 
finite element approximation of this problem and the stabilization of the discretization 
by using the SUPG method. In Section \ref{Sec:dwr} our a posteriori error control 
mechanism based on the DWR method is developed and localized error terms are derived. In 
Section \ref{Sec:PracAsp} some implementational issues are addressed. Finally, in Section 
\ref{Sec:NumStud} the results of numerical computations are presented in order to 
illustrate the feasibility, potential and benefit of the proposed approach. Further, a 
careful comparison with reference values of the literature is given for a benchmark 
problem.

\section{Problem formulation and stabilized discretization}
\label{Sec:pro}

In this work we study the linear convection-diffusion-reaction problem 
\begin{equation}
\label{Eq:pro}
\begin{array}{r@{\;}c@{\;}l@{\hspace*{2ex}}l}
\partial_t u + \vec b\cdot\nabla u - \DIV\left(\eps \nabla u\right) + \alpha u &= & f 
 & \mbox{in } \Omega \times (0, T]\,,\\[1ex]
u(\vec x,t) &= & 0  & \mbox{on } \partial\Omega\times(0,T]\,,\\[1ex]
u(\vec x,0) &= & u_0 & \mbox{in } \Omega\,.
\end{array}
\end{equation}
We assume that $\Omega\subset \R^d$, with $d=2$ or $d=3$, is a polygonal or 
polyhedral bounded domain and let $I:=(0,T]$. 
For simplicity, problem \eqref{Eq:pro} is equipped with homogeneous boundary 
conditions. Problem \eqref{Eq:pro} is considered as a prototype model for more 
sophistcated equations of practical interest, for instance, for the  
Navier--Stokes equations of incompressible viscous flow. For an application of 
our approach to semilinear problems with nonlinear reactive terms we refer 
to \cite{S14}.  

To ensure the well-posedness of problem \eqref{Eq:pro} we assume that 
$\alpha\in\R^+_0$, $\eps \in\Ls^{\infty}(\Omega)$, 
$\vec b\in\vec\Hs^1(\Omega)\cap\vec\Ls^{\infty}(\Omega)$ with 
$\eps(\boldsymbol{x}) \geq \eps_0 > 0$ and 
$\left(\DIV \vec b \right)(\boldsymbol{x}) = 0$ almost everywhere in $\Omega$. 
Further, we let $f\in\Ls^2\left(0,T; \Hs^{-1}_0(\Omega)\right)$ and 
$u_0\in\Ls^2(\Omega)$ with $\Hs^{-1}_0(\Omega)$ denoting the dual space 
of $\Hs^1_0(\Omega)$. Then the existence and uniqueness of a weak solution 
\begin{equation}
\label{Def:X}
u\in \mathcal X := \left\{ v \in \Ls^2\left(0,T;\Hs^1_0(\Omega)\right) \mid \partial_t v 
\in \Ls^2\left(0,T; \Hs^{-1}(\Omega)\right)\right\}
\end{equation}
of problem \eqref{Eq:pro}, satisfying $u(0)=u_0$ and 
\begin{equation}
\label{Eq:WPpw}
\langle \partial_t u, \varphi \rangle + \langle \vec b \cdot \nabla u, \varphi 
\rangle  + \langle \varepsilon \nabla u , \nabla \varphi \rangle  + 
\langle \alpha u , \varphi \rangle  = \langle f, \varphi \rangle 
\end{equation}
for all $\varphi \in  \mathcal H^1_0(\Omega)$ and almost every $t \in (0,T)$, is 
ensured \cite{E10,QV08}. 

The Dual Weighted Residual approach (for short DWR method) is based on a variational 
space-time discretization of problem \eqref{Eq:WP} and of a corresponding adjoint (or 
dual) problem. In the analysis that is given below we need discontinuous and continuous 
variational discretizations of the time variable. The discretization and the 
application of the DWR method is done for variational space-time approximations 
with piecewise polynomials of arbitrary order in space and time. In our numerical 
calculations (cf.\ Section \ref{Sec:NumStud}) we restrict ourselves to applying lowest 
order members of these families of time discretization schemes to the primal and dual 
problem. In this work we aim to demonstrate and analyze the feasibility of our approach to 
convection-dominated transport. For higher order variational discretizations of the time 
variable the solution of the arising algebraic systems of equations becomes much more 
involved. For their application in non-adaptive computations we refer 
to, e.g., \cite{BK15,K15,BK15_2}. 

For the discretization in time we divide the time interval $I$ into not necessarily 
equidistant subintervals $I_m:=(t_{m-1}, t_m]$ , with $m=1,\ldots,M$, where $0 = t_0 < 
t_1 < \ldots < t_{M-1} < t_M = T$ with step size $k_m = t_m - t_{m-1}$ and $k = 
\max\limits_m k_m$. We put 
\begin{equation}
\label{Def:Y}
\mathcal Y := \left\{\Ls^2\left(0,T;\Hs^1_0(\Omega)\right)\mid v_{|I_m} \in 
C(\overline I_m;\Hs^1_0(\Omega)) \right\}\,.
\end{equation}
Here, the notation $v_{|I_m} \in C(\overline I_m;\Hs^1_0(\Omega))$ means that $v_{|I_m}$ 
posses a continuous extension to the closure $\overline I_m$ of $I_m$. The unique 
solution $u\in \mathcal X$ of problem \eqref{Eq:WPpw} then satisfies the variational 
space-time problem: \emph{Find $u \in \mathcal X$ 
such that}
\begin{equation}
\label{Eq:WP}
\langle\langle \partial_t u, \varphi\rangle\rangle + a(u)(\varphi) + \langle 
u(0),\varphi(0) \rangle = F(\varphi) + \langle u_0,\varphi(0) \rangle 
\end{equation}
\emph{for all $\varphi \in  \mathcal Y$.}  

In \eqref{Eq:WP} we use the notation 
\begin{equation}
\label{Def:as}
a(v)(\varphi) := \langle\langle \vec b\cdot\nabla v, 
\varphi \rangle\rangle + \langle\langle \eps\nabla v, \nabla\varphi\rangle\rangle + 
\langle\langle\alpha v,\varphi\rangle\rangle
\end{equation}
for $v, \varphi \in \mathcal Y$ and 
\begin{equation}
\label{def:F}
F(\varphi):= \langle\langle f,\varphi\rangle\rangle 
\end{equation}
for $\varphi \in \mathcal Y$, where  
\[
\langle\langle v,w\rangle\rangle:=\int_0^T \langle v, w\rangle \ud t
\]
is the inner product of $\mathcal L^2(0,T;\mathcal L^2(\Omega))$ and $\langle \cdot, 
\cdot \rangle$ is the 
inner product of $\mathcal L^2(\Omega)$. In \eqref{Eq:WP} the initial condition is 
imposed in a 
weak form. We note that $\mathcal Y$ is a dense subspace of 
$\Ls^2\left(0,T;\Hs^1_0(\Omega)\right)$.

Next, we introduce the time discrete function spaces
\begin{align}
\label{eq:08a}
\mathcal X_k^r & :=\left\{v_k\in\Ls^2\left(0,T; \Hs^{1}_0(\Omega)\right)\mid {v_k}_{\vert 
I_m}\in\Ps_r\left(I_m; \Hs^{1}_0(\Omega)\right), \; v_k(0) \in 
\mathcal{L}^2(\Omega)\right\}\,,\\[1ex]
\label{eq:08b}
\overline{\mathcal X}_k^{\,r} & :=\left\{v_k\in C\left([0,T]; \mathcal L^2 
(\Omega)\right)\mid {v_k}_{\vert 
I_m}\in\Ps_r\left(
\overline I_m
; \Hs^{1}_0(\Omega)\right)\right\}\,,
\end{align}
where $\Ps_r(I_m;\Hs^{1}_0(\Omega))$ denotes the space of all polynomials in time up to 
degree $r\geq 0$ on $I_m$ with values in $\Hs^{1}_0(\Omega)$. For some function $v_k\in 
\mathcal X_k^r$  we define the limits $v_{k,m}^\pm$ from above and below of $v_{k}$ at 
$t_m$ as well as their jump at $t_m$ by 
\[
v_{k,m}^\pm := \lim_{s\rightarrow 0} v_k(t_m \pm s) \,, \qquad [v_k]_m := v_{k,m}^+ 
-v_{k,m}^- \,.
\]
For the temporal discretization of the primal problem \eqref{Eq:WP} we use the 
discontinuous Galerkin method (for short dG($r$)); cf.\ \cite{T06}. The time-discrete 
variational approximation of problem \eqref{Eq:WP} then reads as 
follows: \emph{Find $u_k\in\mathcal X_k^r$ such that}
\begin{align}
\label{eq:09}
A(u_k)(\varphi_k) + \left\langle u_{k,0}^{+} , \varphi_{k,0}^{+} \right\rangle&= F(\varphi_k) 
+ \left\langle u_0, \varphi_{k,0}^{+}\right\rangle
\end{align}
\emph{for all $\varphi_k\in\mathcal X_k^r$.}

In \eqref{eq:09} we use the notation 
\begin{equation}
\label{def:A}
\begin{split}
A(v_k)(\varphi_k) := & \sum_{m=1}^{M}  \int_{I_m} \langle \partial_t v_k, 
\varphi_k\rangle \ud t + a(v_k)(\varphi_k) + \sum\limits_{m=2}^{M} \left\langle 
\left[v_k\right]_{m-1},\varphi_{k,m-1}^{+}\right\rangle
\end{split}
\end{equation}
for $v_k, \varphi_k \in \mathcal X_k^r$. We note that the initial condition is 
incorporated into the variational problem. For the derivation and an analysis of the 
dG($r$) semidiscretization of abstract evolution problems in Hilbert spaces we refer to 
\cite{T06}. The dG($r$) method is nonconforming, since by an embedding result (cf.\ 
\cite{DL92}) it holds that $\mathcal X_k^r \not\subset \mathcal X$.

Next, we describe the Galerkin finite element approximation in space of the semidiscrete 
problem \eqref{eq:09}. To this end, we use two- or three-dimensional shape- and 
contact-regular meshes \cite{C78}. By $\mathcal{T}_h = \{K\}$ we denote a conforming 
decomposition of the domain $\Omega$ into triangles in two space dimensions or tetrahedra 
in three space dimensions. Quadrilateral and hexahedral elements can be applied in 
the same way by means of the standard modifications. On $\mathcal{T}_h$ we define the 
function space $\mathcal{V}_h^p\subset \Hs_0^1(\Omega)$ by 
\[
\mathcal{V}_h^p := \{v\in\Hs_0^1(\Omega)\cap \mathcal{C}(\bar{\Omega})\mid v_{|K} \in 
\mathcal{P}_p(K)\; \forall K\in\mathcal{T}_h\}\,,
\]
with $\mathcal{P}_p(K)$ denoting the function space of polynomials of degree at most $p$ 
on $K$. By replacing $\Hs_0^1(\Omega)$ in the definition of the semidiscrete function 
spaces $\mathcal{X}_k^r$ and $\mathcal{\overline X}_{\, k}^r$ in \eref{eq:08a} and 
\eref{eq:08b}, 
respectively, by $\Vs_h^p$, we obtain the fully discrete 
function space
\begin{align}
\label{Def:Xkhrp}
\mathcal X_{kh}^{r,p} :=\left\{v_{kh}\in\mathcal X_k^r \mid v_{\vert 
I_m}\in\Ps_r(I_m;\Vs_h^{p,m}), \, \mbox{for } m = 1,\ldots, M \,, \; v_{kh}(0) 
\in \mathcal{V}_h^{p}\right\} 
\end{align}
with $\mathcal X_{kh}^{r,p}\subset \mathcal{X}_k^r$ . We note that the spatial finite 
element space $\Vs_h^{p,m}$ is allowed to be different on 
all intervals $I_m$ which is natural in the context of a discontinuous Galerkin 
approximation of the time variable and allows dynamic mesh changes in time. Throughout 
the time steps $k_m$ are kept constant in space. The fully discrete discontinuous in 
time scheme that is studied below then reads as follows: \emph{Find $u_{kh}\in\mathcal 
X_{kh}^{r,p}$, such that}
\begin{align}
\label{Eq:dGcG}
A(u_{kh})(\varphi_{kh})+ \left\langle u_{kh,0}^{+} , \varphi_{kh,0}^{+} \right\rangle&= 
F(\varphi_{kh}) + \left\langle u_0, \varphi_{kh,0}^{+} \right\rangle
\end{align} 
\emph{for all $\varphi_{kh} \in \mathcal X_{kh}^{r,p}$ with $A(\cdot)(\cdot)$ and 
$F(\cdot)$ being defined in \eqref{def:A} and \eqref{def:F}, respectively.} 

In the DWR approach a continuous Galerkin approximation of the time variable is also 
needed. This type of discretization is applied below to the adjoint problem 
of \eqref{Eq:WP}. Here we introduce the continuous Galerkin approximation of the time 
variable and the resulting fully discrete finite element method for the primal problem 
\eqref{Eq:WP} in order to illustrate its definition. The formulation of the continuous 
Galerkin approximation on dynamically changing meshes is more involved since the global 
continuity of functions in the trial space has to be ensured. Let $\{\tau_0,\ldots , 
\tau_r\}$ be a basis of $\P_r(I_m;\R)$ that satisfies the conditions 
\[
\tau_0(t_{m-1})=1\,, \quad \tau_0(t_m)=0\,, \quad \tau_i(t_{m-1})=0 \,, \;\; i=1,\ldots 
,r\,.
\]
Then we define
\[
\mathcal X_{kh}^{r,p,m} = \mathrm{span}\left\{\tau_iv_i \mid v_0 \in 
\mathcal{V}_h^{p,m-1}\,, \; v_i \in \mathcal V_h^{p,m}\,, \; i=1,\ldots , r \right\}
\]
and 
\[
\overline{\mathcal{X}}_{kh}^{\,r,p} := \left\{v_{kh} \in 
C(\overline{I};\mathcal L^2(\Omega)) 
\mid v_{kh}{}_{|I_m} \in \mathcal X_{kh}^{r,p,m} \right\} \subset 
\overline{\mathcal{X}}_k^{\,r}\,.
\]
We note that this definition of the trial space $\overline{\mathcal{X}}_{kh}^{\,r,p}$ 
ensures the continuity of its functions. This is due to the fact that the vanishing 
spatial degrees of freedom in $\mathcal V_h^{p,m-1}$ are coupled only with the temporal 
basis function $\tau_0$ that vanishes in the right end endpoint $t_m$ of $I_m$. The fully 
discrete continuous in time scheme that is applied below then reads as follows:  
\emph{Find 
$u_{kh}\in \overline{\mathcal X}_{kh}^{\, r,p}$, such that}
\begin{equation}
\label{Eq:cGPcG}
\begin{split}
\langle\langle \partial_t u_{kh}, \varphi_{kh} \rangle\rangle + a(u_{kh})(\varphi_{kh}) + 
\langle u_{kh}(0), \varphi_{kh,0}^{+} \rangle =  \langle\langle f,\varphi_{kh} 
\rangle\rangle 
+ \langle u_0, \varphi_{kh,0}^{+} \rangle
\end{split}
\end{equation} 
\emph{for all $\varphi_{kh} \in \mathcal X_{kh}^{r-1,p}$.}

In \eqref{Eq:cGPcG} the initial condition is imposed in a weak form. This scheme 
belongs to the class of Petrov--Galerkin methods since 
the spaces for the trial and test functions differ.    

In this work we focus on convection-dominated problems with small diffussion parameter $0 
< \varepsilon_0 \ll 1$. Then the finite element approximation needs to be stabilized in 
order to reduce spurious and non-physical oscillations of the discrete solution arising 
close to layers. Here, we apply the streamline upwind Petrov--Galerkin method (for short 
SUPG); cf.\ \cite{RST08,JS08,BS12}. The stabilized variant of the fully discrete scheme 
\eqref{Eq:dGcG} then reads as follows: \emph{Find $u_{kh}\in\mathcal X_{kh}^{r,p}$ such 
that}
\begin{align}
\label{Eq:dGcGs}
A_S(u_{kh})(\varphi_{kh}) + \left\langle u_{kh,0}^{+} , \varphi_{kh,0}^{+} \right\rangle &= 
F(\varphi_{kh}) + \left\langle u_0, \varphi_{kh,0}^{+} \right\rangle
\end{align}
\emph{for all $\varphi_{kh}\in\mathcal X_{kh}^{r,p}$.}

In \eqref{Eq:dGcGs} we put
\begin{equation*}
A_S(u_{kh})(\varphi_{kh}) := A(u_{kh})(\varphi_{kh}) + S(u_{kh})(\varphi_{kh})
\end{equation*}
with
\begin{equation}
 \label{Eq:DefStab}
\begin{split}
S(u_{kh})(\varphi_{kh}) &:= \sum_{m=1}^M\int_{I_m} \sum\limits_{K\in\T_h}\delta_K\langle 
R(u_{kh}), \vec  b\cdot\nabla\varphi_{kh} \rangle_K \ud t\\
&\relph{=} + \sum\limits_{m = 2}^{M}\sum\limits_{K\in\T_h}\delta_K
\left\langle\left[u_{kh}\right]_{m-1}, \vec 
b\cdot\nabla\varphi_{kh,m-1}^+\right\rangle_{K}\\[1ex]
&\relph{=} +\sum\limits_{K\in\T_h}\delta_K \left\langle u_{kh,0}^{+} - u_0,\vec 
b\cdot\nabla\varphi_{kh,0}^{+} \right\rangle_{K}\,,\\[2ex]
R(u_{kh}) &:= \partial_t u_{kh} + \vec b\cdot\nabla u_{kh} - \DIV\left(\eps\nabla 
u_{kh}\right)  + \alpha u_{kh} - f\,.
\end{split}\end{equation}
for $\{u_{kh}, \varphi_{kh}\}\in \mathcal X_{kh}^{r,p}\times \mathcal X_{kh}^{r,p}$. In 
Eq.\ \eqref{Eq:DefStab} we denote by $\langle \cdot, \cdot \rangle_{K}$ 
the inner product  of the space $\mathcal L^2(K)$. The proper choice of the stabilization 
parameter $\delta_K$ is an important issue in the application of the SUPG approach; cf.\ 
\cite{JK13} and the discussion  therein. As proposed by our analysis of stabilized finite 
element methods in \cite{BS12} we choose 
\[
\delta_K \sim \min\left\{\dfrac{h_K}{p \| \vec b\|_{\vec{\mathcal 
L}^\infty(K)}};\dfrac{h_K^2}{p^4 \, \varepsilon};\dfrac{1}{
k_m+\alpha};\dfrac{k_m+\alpha}{\alpha^2}\right\}\,.
\]
In \eqref{Eq:DefStab} and from now on we assume for brevity that the diffusion 
coefficient $\varepsilon(x)$ equals a constant $\varepsilon$. Otherwise the  
additional projection operator that is used in \cite{BS12} has to be applied to 
the diffusive term of the residual in \eqref{Eq:DefStab}. The SUPG stabilized 
form of the continuous in time scheme \eqref{Eq:cGPcG} is obtained along the 
same lines.

\begin{rem}
In this work we restrict ourselves to linear problems in the nonstationary case. This 
is sufficient to study and illustrate our dual weighted residual approach for stabilized 
finite element approximations of convection-dominated problems. A further stabilization in 
crosswind direction may be obtained by using an additional shock-capturing stabilization 
technique; cf., e.g.\ \cite{BS12,JS08,JK07}. However, the most 
efficient family of this type of additional stabilization is based on adding additional 
nonlinear terms. In the case of linear problems the latter methods then increase the 
complexity of solving the arising algebraic system of equations significantly. For this 
reason an additional shock-capturing stabilization is not studied here for the 
nonstationary problem \eqref{Eq:pro}. However, some of our numerical studies that are 
presented in Section \ref{Sec:NumStud} are done for steady nonlinear problems. 
Restricting ourselves to steady problems in the nonlinear case is sufficient to 
demonstrate the feasibilty of our method also to nonlinear equations and helps to 
separate characteristic features that are related to the discretization in space. 
\end{rem}

For the steady counterpart of problem \eqref{Eq:pro}, 
\begin{equation}
\label{Eq:pros}
\alpha u + \vec b\cdot\nabla u -\nabla \cdot \left(\eps \nabla u\right) + r(u) = f 
\quad \mbox{in } \Omega, \qquad u = 0 \quad \mbox{on } \partial\Omega\,,
\end{equation}
with the above assumptions about the data and, further, supposing that (cf.\ \cite{BS12}) 
\begin{equation*}
\label{Eq:AssNL}
 r\in C^1(\R)\,,\quad r(0) = 0\,, \quad r'(s)\geq r_0\geq 0 \quad \mbox{for } s\geq 0\,, 
\mbox{ }s\in \R \,,
\end{equation*}
we consider using SUPG and additional shock-capturing stabilization (cf.\ 
\cite{BS12,JK07}) such that the fully discrete problem reads as: \emph{Find $u_h\in 
\mathcal V_h^p$ such that} 
\begin{align}
\label{eq:08}
A_{SC}(u_h)(\varphi_h)&= F(\varphi_h) 
\end{align}
\emph{for all $\varphi_h\in\Vs_h^p$ with}
\begin{equation}
\label{Eq:DefASC}
\begin{split}
A_{SC}(u_h)(\varphi_h) &:= A(u_h)(\varphi_h) + S(u_h)(\varphi_h) + 
S_C(u_h)(\varphi_h)\,,\\[1ex]
A(u_h)(\varphi_h) & = \langle \alpha u_h,\varphi_h\rangle + \langle 
\vec b\cdot\nabla u_h, \varphi_h\rangle + \langle \eps\nabla u_h, 
\nabla\varphi_h\rangle + \langle r(u_h),\varphi_h \rangle  \\[1ex]
S(u_h)(\varphi_h) &:= \sum\limits_{K\in\T_h} \delta_K \langle R(u_h), \vec 
b\cdot\nabla\varphi_h\rangle_{K}\,,\\[1ex] 
R(u_h) &:= \alpha u_h + \vec b\cdot\nabla u_h -\nabla\cdot(\Pi_K \eps\nabla u_h)  + 
r(u_h) - 
f\,,\\[1ex] 
S_C(u_h)(\varphi_h) &:= \sum\limits_{K\in\T_h}\langle \tau_K(u_h) \vec D \nabla 
u_h,\nabla\varphi_h \rangle_{K}\,.\\[1ex]
F(\varphi_h) &:= \langle f,\varphi_h\rangle 
\end{split}
\end{equation}
\emph{and the stabilization parameter}
\begin{equation}
\label{Eq:DefSCPar}
\begin{split}
\delta_K  & \sim \min\left\{\dfrac{h_K}{p \| \vec b\|_{\vec{
\mathcal L}^\infty(K)}};\dfrac{h_K^2}{p^4\mu_{\mathrm{inv}}^2 
\|\varepsilon\|_{\mathcal L^\infty(K)}};\dfrac{1}{\alpha+r_0};\dfrac{\alpha+r_0}{L_r^2}
\right\} \ , ,
\\[1ex]
\tau_K(u_h) &:= l_K(u_h)\hat{R}_K(u_h) = \frac{l_K(u_h)\Vert R(u_h)\Vert_{\Ls^2(K)}}{| 
u_h|_{\Hs^1(K)} + \kappa_K}\,,\\[1ex]
l_K(u_h) &:= l_0 h_K \max\left\{0, \beta - \frac{2\Vert\eps\Vert_{\Ls^{\infty}(K)}}{h_K 
\hat{R}_K(u_h)} \right\}\,,\quad 
\vec D :=\left\{\begin{array}{@{}ll} \vec I - \dfrac{\vec b\otimes\vec b}{\vert \vec 
b\vert^2}\,, & \vec b \neq \vec 0\,, \\ \vec 0\,, & \vec b = \vec 0\,.\end{array}\right.
\end{split}
\end{equation}

In the set of equations \eqref{Eq:DefSCPar} we denote by $\|\cdot \|_{\mathcal L^\infty  
(K)}$ and $\|\cdot \|_{\mathcal H^1 (K)}$ the usual norms associated with the function 
spaces on the element $K$. For further details regarding the definition and choice of 
the 
stabilization parameter in \eqref{Eq:DefSCPar} we refer to \cite{BS12}.

\section{A dual weighted residual approach for stabilized finite element methods} 
\label{Sec:dwr}

Here we develop our application of the Dual Weighted Residual (for short DWR) method 
(cf.\ \cite{BR03}) to the stabilized finite element 
approximation \eqref{Eq:dGcGs} of problem \eqref{Eq:pro}. The DWR approach aims at an 
error control for an arbitrary quantity of physical interest. This is in contrast to 
standard a posteriori error estimates that typically provide computable upper (and lower) 
bounds in terms of numerically available quantities for the numerical approximation 
errors measured in standard norm, for instance in the natural norm of the discretization 
for that an a priori error analysis is available. The capability of providing an 
error control mechanism for physically relevant quantities offers large potential of the 
DWR approach in engineering sciences. The DWR method is based on duality techniques and an 
additional nonstationary adjoint problem has to be solved which includes the primal 
solution as coefficient. Thus, in each adaptation step of an adaptive algorithm 
numerical approximations to the solution of the primal and dual problem need to be 
computed for the whole time period such that the simulations become numerically 
expensive. Several techniques were proposed to reduce the computational costs for 
determining the approximate dual solution. For this topic we refer to the discussion in  
Section \ref{Sec:Intro}. 

The characteristic feature of most of the existing a posteriori error analyses for 
con\-vection-dominated problems is their non-robustness with respect to the small perturbation 
parameter which then leads to adaptive meshes that are not satisfactory yet. On the 
other hand the DWR method yields an exact representation of the discretization error in 
the target quantity. This observation is the key point of our application of the DWR 
concept. The latter representation depends on the exact dual solution that has 
still to be approximated. For this we use higher order techniques which is in contrast to 
other works of the literatur \cite{BR03}. Thereby we aim at a reduction of approximation 
errors in the sensitive regions of convection-dominated problems with sharp layers and 
strong gradients where approximations and interpolations are highly delicious and 
strongly defective. This approach increases the computational costs for solving the 
adjoint problem, but on the other hand it improves the approximation quality of the 
weights in the a posteriori error control mechanism and, thereby, the effectivity of the 
adaptation process. The proper choice of the weights is considered to be an important 
step in the application of the DWR method to stabilized approximations 
of convection-dominated problems and to deserve careful attention.  Our numerical 
computations (cf.\ Section \ref{Sec:NumStud}) will illustrate the impact of the 
approximation of the dual solution on the approximation quality in the target quantity. 

The DWR approach aims to control the error with respect to some output functional 
$\mathcal J(\cdot)$. This requires a respresentation of an estimate of the difference 
$\mathcal J(u)-\mathcal J(u_{kh})$. Here, 
$\mathcal J(u)$ is  the user-chosen target quantity of physical 
interest. We suppose that the functional $\mathcal J(\cdot)$ is defined on the space $\mathcal Y$ 
introduced in \eqref{Def:Y}, i.e.\ $\mathcal J:\mathcal Y \mapsto \R$. Further, we 
assume that the functional $\mathcal J$ is Fr\'{e}chet differentiable, i.e.\ $\mathcal J'(y) \in 
\mathcal Y'$ for $y\in \mathcal Y$. Moreover, we assume that the directional derivative 
of $\mathcal J$ admits an $\mathcal L^2$ representation such that for any $v \in \mathcal Y$ there 
exists 
some function $j(v) \in\mathcal L^2\left(0,T;\mathcal L^2(\Omega)\right)$ such that 
\begin{equation}
\label{Def:AssumJ}
\mathcal J'(v)(\varphi) = \langle\langle j(v),\varphi\rangle\rangle
\end{equation}
is satisfied for all $\varphi \in \mathcal Y$. If the target functional $\mathcal J$ is 
less regular, involving for example spatial or temporal point-values, then the theory 
developed below can no longer be  applied directly. In this case a regularization of the 
functional may be used to overcome  the lack of regularity. However, the regularization is 
usually only necessary in the  development of the formal framework. On the discrete level 
and in practical computations  the abstract theory often performs successfully even for 
less regular output functionals (cf.\ Section \ref{Sec:NumStud}). 
    
For the derivation of an a posteriori error representation for 
$\mathcal J(u) - \mathcal J(u_{kh})$ we employ the Euler--Lagrange method of 
constrained optimization. We define the Lagrangian functional 
$\mathcal L : \mathcal X \times \mathcal Y \mapsto \R$ by   
\begin{equation}
\label{Eq:EL}
\mathcal{L}(u,z):= \mathcal{J}(u) + F(z) - \langle\langle \partial_t u,z \rangle\rangle 
- a(u)(z) - \langle u(0)-u_0,z(0)\rangle 
\end{equation}
with the target quantity $\mathcal J(\cdot)$ and the forms $a(\cdot)(\cdot)$ and 
$F(\cdot)$ being defined in \eref{Def:as} and \eqref{def:F}, respectively. A stationary 
point $\{u,z\}$ of $\mathcal{L}(\cdot,\cdot)$ on $\mathcal X \times \mathcal Y$ 
is determined by 
\begin{equation}
\label{Eq:EL0}
\mathcal L'(u,z)(\psi,\varphi) = 0 \qquad \text{for all}\; \{\psi,\varphi\}\in 
\mathcal X\times \mathcal Y\,,
\end{equation}
or equivalently by the system of equations 
\begin{align}
\label{Eq:EL1}
\langle\langle \partial_t \psi, z \rangle\rangle 
+ a(\psi)(z) + \langle \psi(0), z(0)\rangle & = \mathcal J'(u)(\psi) \qquad \text{for all}\; 
\psi \in \mathcal X\,,\\[1ex]
\label{Eq:EL2}
\langle\langle \partial_t u,\varphi \rangle\rangle 
+ a(u)(\varphi) + \langle u(0),\varphi(0)\rangle & = F(\varphi) + \langle 
u_0,\varphi(0)\rangle \qquad \text{for all}\; \varphi \in \mathcal Y\,.  
\end{align}
The second of these equations, the $z$-component of the stationarity condition, is just 
the given primal problem \eqref{Eq:WP}. 
Equation \eqref{Eq:EL1}, the $u$-component of the stationarity  condition, is called the 
dual or adjoint equation. In particular, the solution $z\in \mathcal Y$ of the adjoint 
problem \eqref{Eq:EL1} can be recovered as the solution of the following variational 
problem: \emph{Find $z\in \mathcal X$ with} 
\begin{equation}
\label{Eq:DP}
\begin{split}
- \langle \langle \partial_t z , \psi\rangle\rangle  - \langle \langle \vec b \cdot 
\nabla z , \psi\rangle \rangle & + \langle \langle \varepsilon \nabla z , 
\nabla \psi\rangle \rangle\\[1ex] 
& + \langle\langle \alpha z , \psi \rangle \rangle + \langle z(T),\psi(T) \rangle = 
\mathcal J'(u)(\psi)
\end{split}
\end{equation}
\emph{for all $\psi \in \mathcal Y$}.

Under the hypothesis \eqref{Def:AssumJ} the dual problem \eqref{Eq:DP} has the structure 
of the primal problem \eqref{Eq:WP} 
but running backward in time. The existence and uniqueness of a solution $z\in \mathcal 
X$ of problem \eqref{Eq:DP} is thus ensured by the same setting and arguments as used 
for the primal problem \eqref{Eq:WP}. For a right-hand side term \eqref{Def:AssumJ} 
and appropriate assumptions about the boundary $\partial \Omega$ of $\Omega$ the 
continuity constraint in the definition \eqref{Def:Y} of $\mathcal Y$ holds (cf.\ 
\cite{E10} and \cite{QV08} for the Sobolev embedding results) such that $z\in \mathcal Y$ 
is ensured. To see that the solution of \eqref{Eq:DP} in fact satisfies the variational 
problem \eqref{Eq:EL1}, we use integration 
by parts with respect to the time variable to find that  
\begin{equation}
\label{Eq:DP_1}
- \langle \langle \partial_t z , \psi\rangle\rangle =  \langle z(0) ,\psi(0)\rangle - 
\langle z(T) ,\psi (T) \rangle +  \langle \langle \partial_t \psi , z \rangle\rangle 
\end{equation}
for test functions $\psi \in \mathcal X$. Combining \eqref{Eq:DP} with \eqref{Eq:DP_1} 
and using integration by parts in the convective term yields \eqref{Eq:EL1}. Below, our 
application of the DWR approach is built upon the dual problem \eqref{Eq:DP}. 

\begin{rem}
In the context of our stabilized finite element approximations two different approaches 
of applying the DWR method can be used. The first approach, refered to as the 
\textbf{first stabilize and then dualize} method, is obtained by introducing a discrete 
Lagrangian functional $\widetilde{\mathcal L}$, that is associated with the stabilized 
Galerkin discretization \eqref{Eq:dGcGs}, and defining the discrete solution 
$\{u_{kh},z_{kh}\}\in\mathcal X_{kh}^{r,p} \times \mathcal X_{kh}^{r,p}$ as the 
stationary point of $\widetilde{\mathcal L}$ on $\mathcal X_{kh}^{r,p}\times \mathcal 
X_{kh}^{r,p}$. To find the desired representation of the error 
$\mathcal J(u) - \mathcal J(u_{kh})$, this quantity is represented in terms of 
the error in the discrete Lagrangian functional; 
cf.\ \cite{SV08}. The second approach, refered to as the \textbf{first dualize and then 
stabilize} method, is obtained by discretizing the continuous Euler--Lagrange system 
\eqref{Eq:EL1}, \eqref{Eq:EL2} by means of the proposed stabilized Galerkin discretization 
scheme \eqref{Eq:dGcGs}, i.e., that the discontinuous in time and continuous in space 
finite element method along with the SUPG stabilization in space is applied to the system 
of equations \eqref{Eq:EL1}, \eqref{Eq:EL2}. As it is shown below, the discrete solution 
$\{u_{kh},z_{kh}\}\in\mathcal X_{kh}^{r,p} \times \mathcal X_{kh}^{r,p}$ is then no 
longer a stationary point of the Lagrangian functional, it's just an approximation to 
such point. In this approach the error in the goal quantity 
$\mathcal J(u) - \mathcal J(u_{kh})$ is 
represented in terms of the continuous Lagrangian functional \eqref{Eq:EL}. The difference 
of the either approaches comes through the presence of the stabilization terms in the 
discrete Lagrangian functional. In this work we apply the second approach. In the second 
approach the SUPG stabilization of the discrete dual problem is based on the residual of 
the discrete counterpart of the backward in time problem \eqref{Eq:DP}; cf.\ Eq.\ 
\eqref{Eq:ELd1} below. This seems to be more natural. Moreover, numerical instabilities 
were observed in the literature \cite{B00} for the first strategy of transposing the whole 
stabilized system. For a careful comparison of the either approaches of applying the DWR 
method to stabilized discretization schemes we refer to \cite{S14} where this is done for 
stationary problems. For illustration purposes we sketch both approaches briefly. Then we 
follow the second one. We note that the resulting numerical scheme differ in general 
since dualization (i.e.\ optimization) and stabilization do not commute. 
\end{rem}

\textbf{First Stabilize and Then Dualize}

The discrete Lagrangian functional $\widetilde{\mathcal L}: \mathcal X_{kh}^{r,p}\times 
\mathcal X_{kh}^{r,p}\mapsto \R$ associated with the stabilized Galerkin discretization 
\eqref{Eq:dGcGs} is defined by 
\begin{equation}
\label{Eq:ELd}
\widetilde{\mathcal L}(u_{kh},z_{kh})  = \mathcal J(u_{kh}) + F(z_{kh}) - 
A_S(u_{kh})(z_{kh}) - \left\langle u_{kh,0}^{+}-u_0,z_{kh,0}^{+}\right\rangle \,.
\end{equation}
A stationary point $\{u_{kh},z_{kh}\}$ of $\widetilde{\mathcal{L}}(\cdot,\cdot)$ on 
$\mathcal X_{kh}^{r,p}\times \mathcal X_{kh}^{r,p}$ is determined by the equation
\[
\widetilde{\mathcal L}{\,'}(u_{kh},z_{kh})(\psi_{kh},\varphi_{kh}) = 0 \qquad \text{for 
all}\; \{\psi_{kh},\varphi_{kh}\}\in \mathcal X_{kh}^{r,p}\times 
\mathcal X_{kh}^{r,p} \,,
\]
or equivalently by the system of equations 
\begin{align*}
A_S(\psi_{kh})(z_{kh}) + \left\langle \psi_{kh,0}^{+},z_{kh,0}^{+}\right\rangle  & = 
\mathcal J'(u_{kh})(\psi_{kh}) \qquad \text{for all}\; 
\psi_{kh} \in \mathcal X_{kh}^{r,p}\,,\\[1ex]
A_S(u_{kh})(\varphi_{kh}) + \left\langle u_{kh,0}^{+},\varphi_{kh,0}^{+}\right\rangle & = 
F(\varphi_{kh}) + \left\langle u_0,\varphi_{kh,0}^{+} \right\rangle \qquad \text{for all}\; 
\varphi_{kh} \in \mathcal X_{kh}^{r,p}\,.  
\end{align*}

\textbf{First Dualize and Then Stabilize}

We discretize the continuous Euler--Lagrange system \eqref{Eq:EL1}, \eqref{Eq:EL2} 
by the proposed stabilized Galerkin discretization scheme \eqref{Eq:dGcGs}. 
Then the identity \eqref{Eq:EL2} yields the discrete primal problem \eqref{Eq:dGcGs}: 
\textit{Find $u_{kh}\in \mathcal X_{kh}^{r,p}$ such that}
\begin{equation}
\label{Eq:ELd1}
A_S(u_{kh})(\varphi_{kh}) + \left\langle u_{kh,0}^{+} , \varphi_{kh,0}^{+} \right\rangle = 
F(\varphi_{kh}) + \left\langle u_0, \varphi_{kh,0}^{+} \right\rangle  \quad \text{for 
all}\; \varphi_{kh} \in \mathcal X_{kh}^{r,p}\,.
\end{equation}
From the continuous dual problem \eqref{Eq:EL1}, rewritten in the form \eqref{Eq:DP}, 
we find by using the proposed stabilized Galerkin discretization scheme \eqref{Eq:dGcGs} 
the following discrete dual problem: \emph{Find $z_{kh}\in\mathcal X_{kh}^{r,p}$ such 
that}
\begin{align}
\label{Eq:ELd2}
A_{S}^\ast (z_{kh})(\psi_{kh}) + \left\langle z_{kh,T}^- , \psi_{kh,T}^- \right\rangle = 
\mathcal J'(u_{kh})(\psi_{kh})\quad \text{for all}\; \psi_{kh} \in \mathcal 
X_{kh}^{r,p}\,.
\end{align}

Further, we define $z_{kh}(0)=z_{k,0}^+$. In \eqref{Eq:ELd2} we put 
\begin{equation*}
A_{S}^\ast(z_{kh})(\psi_{kh}) := A^\ast (z_{kh})(\psi_{kh}) + S^\ast(z_{kh})(\psi_{kh})
\end{equation*}
with
\begin{equation}
\label{def:As}
\begin{split}
A^\ast(z_{kh})(\psi_{kh}) := & \sum_{m=1}^{M}  - \int_{I_m} \langle \partial_t z_{kh}, 
\psi_{kh}\rangle \ud t - \langle \langle b\cdot \nabla z_{kh},\psi_{kh}\rangle \rangle + 
\langle\langle \varepsilon \nabla z_{kh}, \nabla \psi_{kh}\rangle \rangle \\
& + \langle \langle \alpha z_{kh}, \psi_{kh} \rangle \rangle - \sum\limits_{m=2}^{M} 
\left\langle 
\left[z_{kh}\right]_{m-1},\psi_{kh,m-1}^{-}\right\rangle
\end{split}
\end{equation}
and
\begin{equation*}
\begin{split}
S^\ast(z_{kh})(\psi_{kh}) &:= \sum_{m=1}^M\int_{I_m} 
\sum\limits_{K\in\T_h}\delta_K^\ast \langle 
R^\ast(z_{kh}), -\vec b\cdot\nabla\psi_{kh} \rangle_K \ud t\\[1ex]
&\hspace*{-1.6cm }
- \sum\limits_{m = 2}^{M}\sum\limits_{K\in\T_h}\delta_K^\ast\left\langle\left[z_{kh}\right]_{m-1}, 
-\vec b\cdot\nabla\psi_{kh,m-1}^{-}\right\rangle_{K}+\sum\limits_{K\in\T_h}
\delta_K^\ast 
\left\langle z_{kh,M}^- ,- \vec b\cdot\nabla\psi_{kh,M}^- \right\rangle_{K}\,,\\[2ex]
R^\ast(z_{kh}) &:= -\partial_t z_{kh} -\vec b\cdot\nabla z_{kh} - \DIV\left(\eps\nabla 
z_{kh}\right)  + \alpha z_{kh} - j(u_{kh})\,.
\end{split}\end{equation*}
In the definition of the local residual $R^{\ast}(z_{kh})$ we use the assumption (cf.\ 
\eqref{Def:AssumJ}) that $\mathcal J'(u_{kh})(\cdot)$ admits an $\mathcal L^2$ respresentation 
such 
that $\mathcal J'(u_{kh})(\psi_{kh}) = \langle \langle j(u_{kh}),\psi_{kh}\rangle \rangle$ is 
satisfied for all $\psi_{kh} \in \mathcal{X}_{kh}^{r,p}$ with some function $j(u_{kh}) \in 
\mathcal L^2(0,T;\mathcal L^2(\Omega))$. 

To derive a representation of the error $\mathcal J(u) - \mathcal J(u_{kh})$ we 
need some abstract results. 
For this we need to extend the definition of the Lagrangian functional to 
arguments of $(\mathcal X + \mathcal X_{kh}^{r,p})\times \mathcal Y$. 
In the following we let $ \mathcal L:(\mathcal X+ \mathcal X_{kh}^{r,p})\times \mathcal Y$ 
be defined by 
\begin{equation}
\label{Def:Lext}
\begin{split}
\mathcal L (u,z) := &  \mathcal  J(u)+F(z)- \sum_{m=1}^M \int_{I_m} \langle \partial_t 
u ,z\rangle \ud t - a(u)(z)\\[0.5ex] 
& - \sum\limits_{m=2}^{M} \left\langle 
\left[u\right]_{m-1},z_{m-1}^+\right\rangle - \left\langle u(0) 
-u_0,z(0)\right\rangle\,.
\end{split}
\end{equation}
Then it follows that 
\begin{align}
 \nonumber
& \mathcal L_u (u,z)(\psi) + \mathcal L_z(u,z)(\varphi)\\[1ex] 
\nonumber
& = \mathcal J'(u)(\psi) - \sum_{m=1}^M \int_{I_m} \langle \partial_t \psi,z\rangle \ud t 
- a(\psi)(z) - \sum\limits_{m=2}^{M} \left\langle 
\left[\psi\right]_{m-1},z_{m-1}^+\right\rangle - \left\langle 
\psi(0),z(0)\right\rangle\\[0.5ex]\nonumber
& \quad + F(\varphi) - \sum_{m=1}^M \int_{I_m} \langle \partial_t u ,\varphi\rangle \ud 
t 
- a(u)(\varphi) - \sum\limits_{m=2}^{M} \left\langle 
\left[u\right]_{m-1},\varphi_{m-1}^+ \right\rangle - 
\left\langle u(0)
- u_0,\varphi(0)\right\rangle\\[1.5ex] \nonumber
& = \mathcal J'(u)(\psi) + \sum_{m=1}^M \int_{I_m} \langle \partial_t z, \psi\rangle \ud 
t + \langle \langle \vec b \cdot \nabla z ,\psi \rangle \rangle - \langle \langle 
\varepsilon \nabla z, \nabla \psi \rangle \rangle - \langle \langle \alpha z, \psi \rangle 
\rangle \\[0.5ex] \nonumber
&  \quad + \sum\limits_{m=2}^{M} \left\langle 
\left[z\right]_{m-1},\psi_{m-1}^-\right\rangle 
- \left\langle z(T), \psi(T)\right\rangle + F(\varphi) 
- \sum_{m=1}^M \int_{I_m} \langle \partial_t u ,\varphi\rangle \ud 
t\\[0.5ex]
& \quad - a(u)(\varphi) - \sum\limits_{m=2}^{M} \left\langle 
\left[u\right]_{m-1},\varphi_{m-1}^+\right\rangle - \left\langle 
u(0)-u_0,\varphi(0)\right\rangle
\label{Eq:LExt}
\end{align}
for all $\{\psi,\varphi\}\in (\mathcal X + \mathcal X_{kh}^{r,p})\times \mathcal Y$. 

For the stationary point $\{u,z\}$ of $\mathcal{L}(\cdot,\cdot)$ on $\mathcal X 
\times \mathcal Y$ that is determined by \eqref{Eq:EL0} or \eqref{Eq:EL1}, 
\eqref{Eq:EL2}, respectively, we have that $u,z\in C([0,T];\mathcal L^2(\Omega))$. 
Therefore it 
follows that   
\begin{equation}
\label{Eq:SC}
\mathcal L_u (u,z)(\psi) + \mathcal L_z(u,z)(\varphi) = 0 
\end{equation}
for all $\{\psi,\varphi\}\in \mathcal X\times \mathcal Y$. The discrete 
solution $\{u_{kh},z_{kh}\} \in \mathcal X_{kh}^{r,p}\times \mathcal 
X_{kh}^{r,p}$ then satisfies
\begin{equation}
\begin{split}
\label{Eq:DSSC}
\mathcal L_u (u_{kh},& z_{kh})(\psi_{kh}) +  \mathcal 
L_z(u_{kh},z_{kh})(\varphi_{kh})
 =  S(u_{kh})(\varphi_{kh}) + S^\ast(z_{kh})(\psi_{kh}) 
\end{split}
\end{equation}
for all $\{\psi_{kh},\varphi_{kh}\}\in \mathcal X_{kh}^{r,p}\times \mathcal 
X_{kh}^{r,p}$. For the defect of the discrete solution in the stationarity condition 
\eqref{Eq:DSSC} we use the notation 
\begin{equation}
\mathcal D(x_{kh})(y_{kh}) := S(u_{kh})(\varphi_{kh}) + S^\ast(z_{kh})(\psi_{kh}) 
\end{equation}
with $x_{kh}:= \{u_{kh},z_{kh}\} \in \mathcal X_{kh}^{r,p}\times \mathcal 
X_{kh}^{r,p}$ and $y_{kh}:=\{\psi_{kh},\varphi_{kh}\}\in \mathcal 
X_{kh}^{r,p}\times \mathcal 
X_{kh}^{r,p}$. 

To derive a representation of the error $\mathcal J(u) - \mathcal J(u_{kh})$ 
we need the following abstract theorem that develops the error in terms of the 
Lagrangian functional; cf.\ \cite{S14}.

\begin{theorem}
\label{Thm:L} 
Let $X$ be a function space and $\mathcal L:X \mapsto R$ be a three times differentiable 
functional on $X$. Suppose that $x_c\in X_c$ with some (''continuous'') function space 
$X_c \subset X$ is a stationary point of $\mathcal L$. Suppose that $x_d \in X_d$ with 
some (''discrete'') function space $X_d \subset X$, with not necessarily $X_d \subset 
X_c$, is a Galerkin approximation to $x_c$ being defined by the equation 
\begin{equation}
\label{Def:L0}
\mathcal L'(x_d)(y_d) = \mathcal D(x_d)(y_d) 
\end{equation}
for all $y_d\in X_d$. In addition, suppose that the auxiliary condition 
\begin{equation}
\label{Def:L1}
\mathcal L'(x_c)(x_d) = 0 
\end{equation}
is satisfied. Then there holds the error representation
\begin{equation}
\label{Def:L3}
\mathcal L(x_c) - \mathcal L(x_d) = \frac{1}{2} \mathcal L'(x_d) (x_c- y_d) + \frac{1}{2} 
\mathcal D(x_d)(y_d -x_d) + \mathcal R \,,
\end{equation}
for all $y_d\in X_d$, where the remainder $\mathcal R$ is defined by
\begin{equation}
\label{Def:L4}
\mathcal R = \frac{1}{2}\int_0^1 \mathcal L'''(x_d +s e)(e,e,e) \cdot s \cdot (s-1)\ud s
\end{equation}
with the notation $e:=x_c-x_d$.
\end{theorem}

\begin{proof} In order to keep this work self-contained the proof of Theorem 
\ref{Thm:L} is given in the appendix.
\end{proof}

We note that Theorem \ref{Thm:L} differs from similar theorems that are 
presented in \cite{BR03,BGR10,SV08}, for instance, since in our case the 
discrete solution $\{u_{kh},z_{kh}\}\in\mathcal X_{kh}^{r,p} \times 
\mathcal X_{kh}^{r,p}$ is not a stationary point of a Lagrangian functional but only an 
approximation to such point. In our case the assumption \eqref{Def:L1} is fulfilled by 
means of \eqref{Eq:LExt} along with the definition of the function spaces yielding that 
$\mathcal X_{kh}^{r,p}\subset \mathcal Y$. Theorem \ref{Thm:L} now enables us to derive 
an error representation in terms of the target quantity $\mathcal J(\cdot)$. Here we do 
not separate the error of the temporal and spatial discretization. We study 
directly the error between the continuous and the fully discrete solution which 
is in contrast to the approach in \cite{SV08} for instance.  

For the representation of the error in terms of the target quantity 
$\mathcal J(\cdot)$ we still define the primal residual $\rho(u_{kh})(\cdot)$ 
and the adjoint residual $\rho^\ast(z_{kh})(\cdot)$ by means of 
\begin{align}
\label{Def:PRes}
\rho(u_{kh})(\varphi) & := F(\varphi) - A(u_{kh})(\varphi)- \left\langle 
u_{kh,0}^{+} -u_0,\varphi(0)\right\rangle\\[1ex]
\label{Def:DRes}
\rho^\ast(z_{kh})(\psi) & := \mathcal J'(u_{kh})(\psi) - A^\ast(z_{kh})(\psi)- 
\left\langle 
z_{kh,M}^{-}, \psi(T) \right\rangle
\end{align}
for arbitrary $\varphi \in \mathcal Y$ and $\psi \in \mathcal X + \mathcal X_{kh}^{r,p}$.

\begin{theorem}
\label{Thm:J}
Suppose that $\{u,z\}\in \mathcal{X}\times \mathcal Y$ is a stationary point of 
the Lagrangian functional $\mathcal L$ defined in \eqref{Def:Lext} such that 
\eqref{Eq:SC} is satisfied. Let $\{u_{kh},z_{kh}\}\in\mathcal 
X_{kh}^{r,p} \times \mathcal X_{kh}^{r,p}$ denote its Galerkin approximation being 
defined by \eqref{Eq:ELd1} and \eqref{Eq:ELd2} such that \eqref{Eq:DSSC} is satisfied. 
Then there holds the error representation that 
\begin{equation}
\label{Thm:J1}
\mathcal J(u) - \mathcal J(u_{kh}) = \frac{1}{2}\rho(u_{kh}) (z - \varphi_{kh}) + 
\frac{1}{2}\rho^\ast(z_{kh}) (u-\psi_{kh}) + \mathcal R_{\mathcal S} + 
\mathcal R_{\mathcal J}
\end{equation}
for arbitrary functions $\{\varphi_{kh},\psi_{kh}\}\in \mathcal 
X_{kh}^{r,p}\times \mathcal X_{kh}^{r,p}$, where the 
remainder terms are defined by 
\begin{equation}
\label{Thm:J2}
\mathcal R_{\mathcal S} := \frac{1}{2} S(u_{kh})(\varphi_{kh}+z_{kh}) + \frac{1}{2} 
S^\ast(z_{kh})(\psi_{kh}-u_{kh})
\end{equation}
and
\begin{equation}
\label{Thm:J3}
\mathcal R_{\mathcal J} := \frac{1}{2}\int_0^1 \mathcal J'''(u_{kh}+s\cdot e)(e,e,e) 
\cdot s\cdot (s-1)\ud s
\end{equation}
with $e=u-u_{kh}$.
\end{theorem}

\begin{proof}
Let $x:=\{u,z\}$ with $\{u,z\}\in \mathcal{X}\times \mathcal Y$ be a 
stationary point of $\mathcal L$ in \eqref{Def:Lext} such that \eqref{Eq:SC} is 
satisfied. Let $x_{kh}:=\{u_{kh},z_{kh}\}$ with $\{u_{kh},z_{kh}\}\in\mathcal 
X_{kh}^{r,p} \times \mathcal X_{kh}^{r,p}$ denote the Galerkin approximation of $x$ 
that is defined by \eqref{Eq:ELd1} and \eqref{Eq:ELd2}, respectively. From 
\eqref{Def:Lext} along with \eqref{Eq:SC} and \eqref{Eq:ELd1} we conclude that 
\begin{equation*}
\label{Thm:J4}
\mathcal J(u) - \mathcal J(u_{kh}) = \mathcal L(x) - \mathcal L(x_{kh}) + 
S(u_{kh})(z_{kh})\,.
\end{equation*}
As mentioned above, condition \eqref{Def:L1} is satisfied in our case. By Thm.\ 
\ref{Thm:L} we get that 
\begin{equation}
\label{Thm:J5}
\mathcal J(u) - \mathcal J(u_{kh}) = \frac{1}{2} \mathcal{L}'(x_{kh})(x-y_{kh}) + 
\frac{1}{2} \mathcal D(x_{kh})(y_{kh} -x_{kh}) + S(u_{kh})(z_{kh})+ \mathcal R
\end{equation}
for all $y_{kh}=\{\psi_{kh},\varphi_{kh}\}\in \mathcal X_{kh}^{r,p} \times \mathcal 
X_{kh}^{r,p}$ with the remainder $\mathcal R$ being defined by \eqref{Def:L4}. Recalling 
the definition \eqref{Def:Lext} of $\mathcal L$ yields for the remainder $\mathcal R$ the 
asserted representation \eqref{Thm:J3}.

Next, from \eqref{Eq:LExt} along with the definitions \eqref{def:A} and \eqref{def:As} it 
follows that
\begin{align*}
& \mathcal{L}'(u_{kh},z_{kh})(u-\psi_{kh},z- \varphi_{kh}) \\ 
& = \mathcal{J}'(u_{kh})(u-\psi_{kh}) - A^{\ast}(z_{kh})(u-\psi_{kh}) - \left\langle 
z_{kh,M}^{-}, u(T)-\psi_{kh,M}^- \right\rangle \\[1ex]
& \qquad + F(z-\varphi_{kh}) -A(u_{kh})(z-\varphi_{kh}) - \left\langle 
u_{kh,0}^{+} -u_0,z(0)- \varphi_{kh,0}^{+}\right\rangle  \\[1.5ex]
& = \rho^\ast(z_{kh})(u-\psi_{kh}) + \rho(u_{kh})(z-\varphi_{kh})
\end{align*}
for all $\{\psi_{kh},\varphi_{kh}\}\in \mathcal X_{kh}^{r,p} \times \mathcal 
X_{kh}^{r,p}$. Substituting this identity into \eqref{Thm:J5} yields that
\begin{equation}
\label{Thm:J6}
\begin{split}
\mathcal J(u) - \mathcal J(u_{kh})  & = \frac{1}{2} 
\rho^\ast(z_{kh})(u-\psi_{kh}) + \frac{1}{2}
\rho(u_{kh})(z-\varphi_{kh})\\
& \qquad + \frac{1}{2} \mathcal D(x_{kh})(y_{kh} -x_{kh}) + S(u_{kh})(z_{kh})+ 
\mathcal R_{\mathcal J}\,.
\end{split}
\end{equation}
Finally, we note that 
\begin{align}
\nonumber
& \frac{1}{2} \mathcal D(x_{kh})(y_{kh} -x_{kh}) + 
S(u_{kh})(z_{kh})\\[1ex]\nonumber
& = \frac{1}{2} S(u_{kh})(\varphi_{kh}-z_{kh}) + \frac{1}{2} 
S^\ast(z_{kh})(\psi_{kh}-u_{kh})
+ S(u_{kh})(z_{kh})\\[1.5ex] 
& \begin{gathered} \displaystyle 
= \frac{1}{2} S(u_{kh})(\varphi_{kh}+z_{kh}) + \frac{1}{2} 
S^\ast(z_{kh})(\psi_{kh}-u_{kh})
\,.
\end{gathered}
\label{Thm:J7}
\end{align}
Combining \eqref{Thm:J6} with \eqref{Thm:J7} proves the assertion of the theorem. 
\end{proof}

In the error respresentation \eqref{Thm:J1} the continuous solution $u$ or some 
higher order approximation of $u$ is required for the evaluation of the adjoint 
residual. In the following theorem we show that the adjoint residual coincides 
with the primal residual up to a quadratic remainder. This observation will be 
exploited below to find our final error respresentation in terms of the goal 
quantity $\mathcal J$ and a suitable linearization for its computational evaluation or 
approximation, respectively.

\begin{theorem}
\label{Thm:ResDev}
Suppose that $\{u,z\}\in \mathcal{X}\times \mathcal Y$ is a stationary point of 
the Lagrangian functional $\mathcal L$ defined in \eqref{Def:Lext} such that 
\eqref{Eq:SC} is satisfied. Let $\{u_{kh},z_{kh}\}\in\mathcal 
X_{kh}^{r,p} \times \mathcal X_{kh}^{r,p}$ denote its Galerkin approximation being 
defined by \eqref{Eq:ELd1} and \eqref{Eq:ELd2} such that \eqref{Eq:DSSC} is satisfied. 
Let the primal and adjoint residuals be defined by \eqref{Def:PRes}, \eqref{Def:DRes}. 
Then there holds that 
\begin{equation}
\label{Thm:ResDev1}
\begin{split}
\rho^\ast(z_{kh}) (u-\psi_{kh}) & = \rho(u_{kh}) (z - \varphi_{kh}) + 
S(u_{kh})(\varphi_{kh}-z_{kh}) 
\\[1ex]
& \qquad + S^\ast (z_{kh})(u_{kh}-\psi_{kh}) + \Delta \rho_{\mathcal 
J} 
\end{split}
\end{equation}
for all $\{\psi_{kh},\varphi_{kh}\}\in \mathcal X_{kh}^{r,p}\times \mathcal 
X_{kh}^{r,p}$ with the remainder term
\begin{equation}
\label{Thm:ResDev2}
\Delta \rho_{\mathcal J}:=  - \int_0^1 J''(u_{kh} + s \cdot e)(e,e) \ud s
\end{equation}
with $e:=u - u_{kh}$.
\end{theorem}

\begin{proof}
Let $e:=u - u_{kh}$ and $e^\ast:=z -z_{kh}$ denote the primal and adjoint error, 
respectively. For arbitrary $\psi_{kh}\in \mathcal X_{kh}^{r,p}$ we put
\begin{align*}
k(s) & := \mathcal J'(u_{kh} + s \cdot e)(u-\psi_{kh})- A^\ast(z_{kh}+s 
\cdot e^\ast)(u-\psi_{kh}) \\[1ex]
& \qquad - \left\langle
z_{kh,M}^- + s \cdot e^{\ast,-}_{M}, u(T)-\psi_{kh,M}^- \right\rangle \,.
\end{align*}
We have that 
\begin{align*}
k(1) & :=\mathcal J'(u)(u-\psi_{kh}) - A^\ast(z)(u-\psi_{kh})  - \left\langle
z(T), u(T)-\psi_{kh,M}^- \right\rangle = 0\,.
\end{align*}
From \eqref{Def:DRes} we get that
\begin{align*}
k(0) & = \mathcal J'(u_{kh})(u - \psi_{kh}) - A^\ast(z_{kh})(u-\psi_{kh}) - \left\langle 
z_{kh,M}^-, u(T)-\psi_{kh,M}^- \right\rangle \\[1ex]
& = \rho^\ast (z_{kh}) (u-\psi_{kh})\,. 
\end{align*}
Further, we conclude that 
\begin{align*}
k'(s) & =  \mathcal J''(u_{kh} + s\cdot e)(e,u-\psi_{kh})- 
A^\ast(e^\ast)(u-\psi_{kh})\\[1ex] 
& \qquad - 
\left\langle e^{\ast,-}_{M}, u(T)-\psi_{kh,M}^- \right\rangle \,.
\end{align*}
Using \eqref{Eq:ELd2} and \eqref{Def:DRes} we find that 
\begin{align}
\nonumber
& \rho^\ast (z_{kh}) (u-\psi_{kh})  = \mathcal J'(u_{kh})(u - \psi_{kh}) - 
A^\ast(z_{kh})(u-\psi_{kh})\\[1ex] \nonumber
& \qquad - \left\langle z_{kh,M}^-, u(T)-\psi_{kh,M}^- \right\rangle + S^\ast 
(z_{kh})(\psi_{kh}) - S^\ast (z_{kh})(\psi_{kh}) \\[1ex] \nonumber
& \qquad - \mathcal J'(u_{kh})(u_{kh}) + A^\ast(z_{kh})(u_{kh}) + 
S^\ast(z_{kh})(u_{kh}) + \left\langle z_{kh,M}^-, u_{kh,M}^- 
\right\rangle\\[1.5ex]\nonumber 
& = \;  \rho^\ast (z_{kh}) (u-u_{kh}) +  S^\ast 
(z_{kh})(u_{kh}-\psi_{kh})\\[1ex]\label{Thm:ResDev5}
& \; =  \rho^\ast (z_{kh}) (e) +  S^\ast (z_{kh})(u_{kh}-\psi_{kh})\,.
\end{align}
From \eqref{Thm:ResDev5} along with the theorem of calculus $\int\limits_0^1 k'(s) \ud s= 
k(1)-k(0)$ it follows that 
\begin{align}
\nonumber 
& \rho^\ast (z_{kh})(u-\psi_{kh}) = \rho^\ast (z_{kh})(e)  + S^\ast 
(z_{kh})(u_{kh}-\psi_{kh})\\[1ex]\nonumber
&  \;\; = k(0)-k(1)+ S^\ast (z_{kh})(u_{kh}-\psi_{kh}) \\[1ex] 
\nonumber
& \;\; = \int_0^1 \Big(A^\ast(e^\ast)(e) + 
\left\langle e^{\ast,-}_{M}, e_{M}^- \right\rangle - J''(u_{kh} + s \cdot e)(e,e) 
\Big)\ud s+ S^\ast(z_{kh})(u_{kh}-\psi_{kh})\\[1ex]\label{Thm:ResDev6}
& \;\; = A^\ast(e^\ast)(e) + \left\langle e^{\ast,-}_{M}, e_{M}^- \right\rangle 
+ S^\ast(z_{kh})(u_{kh}-\psi_{kh}) + \Delta \rho_{\mathcal J}\,.
\end{align}
Next, for the first and second of the terms on the right-hand side of \eqref{Thm:ResDev6} 
we get that 
\begin{align}
\nonumber
A^\ast(e^\ast)(e) + \left\langle e^{\ast,-}_{M}, e_{M}^- \right\rangle = & 
\; \sum_{m=1}^{M}  - \int_{I_m} \langle \partial_t e^\ast, e \rangle \ud t - \langle 
\langle b\cdot \nabla e^\ast,e \rangle \rangle + \langle\langle \varepsilon \nabla e^\ast, 
\nabla e \rangle \rangle \\ \nonumber
& + \langle \langle \alpha e^\ast, e \rangle \rangle - \sum\limits_{m=2}^{M} 
\left\langle 
\left[e^\ast\right]_{m-1},e_{m-1}^{-}\right\rangle + \left\langle e^{\ast,-}_{M}, e_{M}^- 
\right\rangle \\[1ex]\nonumber
= & \;  \sum_{m=1}^{M} 
 \int_{I_m} \langle \partial_t e, e^\ast \rangle \ud t + \langle \langle b\cdot 
\nabla e, e^\ast \rangle \rangle + \langle\langle \varepsilon \nabla e, \nabla 
e^\ast \rangle \rangle \\ \nonumber
& + \langle \langle \alpha e, e^\ast \rangle \rangle + \sum\limits_{m=2}^{M} 
\left\langle 
\left[e\right]_{m-1},e_{m-1}^{\ast,+} \right\rangle 
+ \left\langle e_0^{+}, e_{0}^{\ast,+} 
\right\rangle\\[1ex]\nonumber
= & \; F(e^\ast) + \left\langle u_0 ,  e_{0}^{\ast,+} \right\rangle - A(u_{kh})(e^\ast)- 
\left\langle u_{kh,0}^{+} ,  e_{0}^{\ast,+} \right\rangle \\[1ex]\label{Thm:ResDev7}
= & \; \rho(u_{kh})(z-z_{kh}) = \rho(u_{kh})(z-\varphi_{kh}) + 
S(u_{kh})(\varphi_{kh}-z_{kh})
\end{align}
for all $\varphi_{kh}\in \mathcal X_{kh}^{r,p}$. Combining \eqref{Thm:ResDev6} with 
\eqref{Thm:ResDev7} yields that
\[
\begin{split}
\rho^\ast (z_{kh})(u-\psi_{kh}) & = \rho(u_{kh})(z-\varphi_{kh}) + 
S(u_{kh})(\varphi_{kh}-z_{kh})\\[1ex]
& \qquad + S^\ast (z_{kh})(u_{kh}-\psi_{kh}) + \Delta \rho_{\mathcal J}
\end{split}
\]
for all $\{\psi_{kh},\varphi_{kh}\}\in \mathcal X_{kh}^{r,p}\times \mathcal 
X_{kh}^{r,p}$ with $\Delta \rho_{\mathcal J}$ being defined by 
\eqref{Thm:ResDev2}. This proves the assertion of the theorem.
\end{proof}

We summarize the results of the previous two theorems in the following corollary. 

\begin{cor}
\label{Cor:J}
Suppose that $\{u,z\}\in \mathcal{X}\times \mathcal Y$ is a stationary point of 
the Lagrangian functional $\mathcal L$ defined in \eqref{Def:Lext} such that 
\eqref{Eq:SC} is satisfied. Let $\{u_{kh},z_{kh}\}\in\mathcal 
X_{kh}^{r,p} \times \mathcal X_{kh}^{r,p}$ denote its Galerkin approximation being 
defined by \eqref{Eq:ELd1} and \eqref{Eq:ELd2} such that \eqref{Eq:DSSC} is satisfied. 
Then there holds the error representation that 
\begin{equation}
\label{Cor:J1}
\mathcal J(u) - \mathcal J(u_{kh}) = \rho(u_{kh}) (z - \varphi_{kh}) + \mathcal 
R_{\mathcal S} + \mathcal R_{\mathcal J} + \frac{1}{2} \Delta \rho_{\mathcal J}
\end{equation}
for arbitrary functions $\varphi_{kh} \in \mathcal X_{kh}^{r,p}$, where the 
primal residual $\rho(u_{kh})(\cdot)$ is defined by \eqref{Def:PRes}, the remainder term 
$\mathcal R_{\mathcal J}$ is given by \eqref{Thm:J3}, the linearization error $\Delta 
\rho_{\mathcal J}$ is defined by \eqref{Thm:ResDev2} and 
\begin{equation}
\label{Cor:J2}
\mathcal R_{\mathcal S} := S(u_{kh})(\varphi_{kh}) 
\end{equation}
for arbitrary functions $\psi_{kh} \in \mathcal X_{kh}^{r,p}$.
\end{cor}

We note that the notation $\mathcal R_{\mathcal S}$ in \eqref{Cor:J2} is used 
generically and defined differently in different equations of its occurence. In 
the final step of deriving an a posteriori error representation we give a 
localized error approximation that can be used to design an adaptive algorithm. 

\begin{theorem}[Localized error representation]
Suppose that $\{u,z\}\in \mathcal{X}\times \mathcal Y$ is a stationary point of 
the Lagrangian functional $\mathcal L$ defined in \eqref{Def:Lext} such that 
\eqref{Eq:SC} is satisfied. Let $\{u_{kh},z_{kh}\}\in\mathcal 
X_{kh}^{r,p} \times \mathcal X_{kh}^{r,p}$ denote its Galerkin approximation being 
defined by \eqref{Eq:ELd1} and \eqref{Eq:ELd2} such that \eqref{Eq:DSSC} is satisfied. 
Neglecting the higher order error terms in \eqref{Cor:J1}, then there holds as a 
linear approximation the cell-wise error representation
\begin{equation}
\label{Thm:CellErr1}
\begin{split}
\J(u) - \J(u_{kh}) & \doteq \int\limits_0^T \sum\limits_{K\in\T_h} \Big\{\langle 
\mathcal{R}(u_{kh}), z-\varphi_{kh}\rangle_{K} - \delta_K\langle 
\mathcal{R}(u_{kh}), \vec b\cdot\nabla\varphi_{kh}\rangle_{K}\\
&\relph{=}  - \langle \mathcal{E}(u_{kh}), z-\varphi_{kh}\rangle_{\partial K}\Big\} \ud t 
- \left\langle  u_{kh,0}^{+} - u_0, z(t_0) - \varphi_{kh,0}^{+} \right\rangle_{\Omega}\\
&\relph{=} - \sum\limits_{m=2}^{M} \left\langle\left[u_{kh}\right]_{m-1}, z(t_{m-1}) - 
\varphi_{kh,m-1}^+ \right\rangle_{\Omega} \\
&\relph{=} + \sum\limits_{K\in\T_h} \delta_K\left\langle u_{kh,0}^{+} - u_0, 
\vec b\cdot\nabla\varphi_{kh,0}^{+}\right\rangle_{K}\\
&\relph{=} + \sum\limits_{m=2}^{M} \sum\limits_{K\in\T_h} \delta_K 
\left\langle\left[u_{kh}\right]_{m-1}, 
\vec b\cdot\nabla \varphi_{kh, m-1}^+\right\rangle_{K} \,.
\end{split}
\end{equation}
The cell- and edge-wise residuals are defined by 
\begin{align}
\label{eq:42} \mathcal R(u_{kh})_{|K} & := 
f - \partial_t u_{kh} - \vec b\cdot\nabla u_{kh} + \nabla\cdot(\eps\nabla u_{kh})  
-\alpha u_{kh} \,,\\[0.5ex]
\label{eq:43} \mathcal E(u_{kh})_{|\Gamma} & := 
\left\{ \begin{array}{cl} \frac{1}{2}\vec n\cdot[\eps\nabla u_{kh}] 
& \mbox{ if } \Gamma\subset\partial K\backslash\partial\Omega\,, \\[0.5ex] 
0 & \mbox{ if } \Gamma\subset\partial\Omega\,,\\ \end{array}\right.
\end{align}
where $[\nabla u_{kh}]:= 
\nabla u_{kh}{}_{|\Gamma\cap K'}-\nabla u_{kh}{}_{|\Gamma\cap K}$ defines the 
jump of $\nabla u_{kh}$ over the inner edges $\Gamma$ with normal unit vector 
$\vec n$ pointing from $K'$ to $K$. 
\end{theorem}

\begin{proof} 
The assertion directly follows from \eqref{Cor:J1}, \eqref{Cor:J2} by neglecting 
the higher order remainder terms $\mathcal R_{\mathcal J}$ and 
$\Delta \rho_{\mathcal J}$
as well as applying integration by parts on each cell $K\in \mathcal{T}_h$ to 
the diffusive term in the primal residual \eqref{Def:PRes}. 
\end{proof}

Finally, we summarize the result of our application of the DWR approach to the 
stabilized approximation \eqref{eq:08}--\eqref{Eq:DefSCPar} of the nonlinear stationary 
problem \eqref{Eq:pros}. In terms of a \emph{first dualize and then stabilize} philosophy 
analogously to \eqref{Eq:ELd2} we get that 
\begin{equation}
\label{Eq:DWRPros}
\J(u) - \J(u_h) = \rho(u_h)(z-\varphi_h) + \mathcal R_{\St} + \mathcal 
R_{\mathrm{nl}}+\frac{1}{2}\Delta\rho_{\St} +\frac{1}{2}\Delta\rho_{\mathrm{nl}}
\end{equation}
with the primal residual 
\begin{equation}
\label{Eq:PResS}
\rho(u_h)(\varphi) := F(\varphi) - A(u_h)(\varphi) 
\end{equation}
and the remainder terms of the stabilization  
\begin{equation}
\label{Eq:RDst}
\mathcal R_{\St} +\frac{1}{2}\Delta\rho_{\St} = S(u_h)(\varphi_h) + S_C(u_h)(\varphi_h)
\end{equation}
as well as the higher order remainder terms
\begin{align}
\nonumber
\mathcal R_{\mathrm{nl}} & :=\frac{1}{2}\int_0^1 
\Big\{\J'''(u_h+se)(e,e,e) - \left\langle r'''(u_h+se)e^3, z_h+se^\ast 
\right\rangle \\[1ex]
\label{Eq:Rnl}
&\relph{=\frac{1}{2}\int_0^1}- 3 \left\langle r''(u_h+se)e^2, e^\ast \right\rangle_\Omega 
\Big\}\cdot s\cdot (s-1)\ud s\,,\\[2ex]
\label{Eq:Deltanl}
\Delta\rho_{\mathrm{nl}} & = S(u_h)(\varphi_h - z_h) + S_C(u_h)(\varphi_h - z_h) - 
S^\ast(u_h)(\zeta_h-u_h, z_h)\,.
\end{align}
The forms arising in \eqref{Eq:PResS} to \eqref{Eq:Deltanl} are defined in 
\eqref{Eq:DefASC}. We denote by $e:=u-u_h$ and $e^\ast:= z-z_h$ the approximation error 
of the primal and adjoint problem, respectively. For a proof of \eqref{Eq:DWRPros}  we 
refer to \cite{S14}. In \cite{S14}, the \emph{first stabilize and then dualize} approach 
to stationary convection-dominated problems is further presented, investigated 
numerically and compared with the error representation \eqref{Eq:DWRPros}. Finally, 
neglecting the higher order remainder terms $\mathcal R_{\mathrm{nl}}$ 
and $\Delta\rho_{\mathrm{nl}}$ defined in \eqref{Eq:Rnl} and \eqref{Eq:Deltanl}, 
respectively, and using integration by parts we derive from \eqref{Eq:DWRPros} the 
linearized cell-wise error representation 
(cf.\ \cite{S14})
\begin{equation}
\begin{split}
\label{Eq:DWRProsLoc}
\J(u) - \J(u_h) &\doteq \sum\limits_{K\in\T_h} \Big\{\langle \mathcal R(u_h), 
z-\varphi_h\rangle_K - \delta_K\langle\mathcal R(u_h), \vec 
b\cdot\nabla\varphi_h\rangle_K\\
&\relph{=\sum\limits_{K\in\T_h}} + S_{C}(u_h)(\varphi_h) - \langle\mathcal E(u_h), 
z-\varphi_h\rangle_{\partial K}\Big\}\,.
\end{split}
\end{equation}
The cell and edge residuals are defined analogously to \eqref{eq:42}, \eqref{eq:43} by 
\begin{align}
\label{Eq:LocErrS_2} 
\mathcal R(u_h)_{|K} &= f + \nabla\cdot(\eps\nabla u_h) - \vec 
b\cdot\nabla u_h -\alpha u_h - r(u_h)\,,\\[1ex]
\label{Eq:LocErrS_3} 
\mathcal E(u_h)_{|\Gamma} &= \left\{ \begin{array}{@{}cl} 
\frac{1}{2}\vec n\cdot[\eps\nabla u_h]\,, & \mbox{ if } \Gamma\subset\partial 
K\backslash\partial\Omega\,, \\[1ex] 0\,, & \mbox{ if } \Gamma\subset\partial\Omega\,.
\end{array}\right.
\end{align}

\section{Practical Aspects}
\label{Sec:PracAsp}

In this section we discuss some aspects of the practical use of the DWR approach 
presented in Section \ref{Sec:dwr} for the numerical approximation of convection-dominated 
problems. The error representation \eqref{Thm:CellErr1}, written in the 
form 
\begin{equation}
\label{Eq:DefEta}
\J(u) - \J(u_{kh}) \doteq \eta :=\sum\limits_{m=1}^M\sum\limits_{K\in\T_h}\eta_K^m\,,
\end{equation}
depends on the discrete primal solution as well as on the exact dual solution $z$. For 
solving the primal problem \eqref{Eq:pro} we use the discontinuous in time scheme 
\eqref{Eq:dGcGs} and compute a discrete solution $u_{kh}\in \mathcal{X}_{kh}^{r,p}$. As 
mentioned in Section \ref{Sec:Intro} and shown by \eqref{Cor:J2} and 
\eqref{Thm:CellErr1}, respectively, the approximation of the dual solution cannot be done 
in the finite element space of the primal problem since it would result in a vanishing 
primal residual $\rho(u_{kh})(\cdot)$. For the approximation of the dual solution $z$ we 
use the SUPG stabilized counterpart of the continuous in time scheme \eqref{Eq:cGPcG} and 
compute a discrete approximation $z_{kh}\in \overline{\mathcal X}_{kh}^{\,r+1,p+1}$. 

In contrast to many other works of the literature we thus use a higher order approach, 
i.e.\ $z_{kh}\in \overline{\mathcal X}_{kh}^{\,r+1,p+1}$, compared with the primal 
problem for the approximation of the dual solution which leads to higher computational 
costs; cf.\ \cite{BK15,K15} for algorithmic formulations and analyses. In the literature 
the application of a higher order interpolation is often suggested for the 
DWR approach; cf.\ \cite{BR03,BGR10}. For convection-dominated problems such an 
interpolation might be defective and lead to tremendous errors close to sharp layers and 
fronts. Higher order techniques show more stability and reduce spurious oscillations (cf.\ 
\cite{BS12}) which is our motivation for using a higher order approximation of the dual 
solution. For nonlinear problems the additional costs for computing the higher order 
approximation of the dual solution are moderate, since the adjoint problem is always a 
linear one and, thereby, does not require nonlinear (e.g.\ Newton) iterations for solving 
the discrete problem whereas such iterations become necessary in the nonlinear case for 
the primal problem.     

In order to define the localized error contributions $\eta_K^m$ in \eqref{Eq:DefEta} we 
consider a hierarchy of sequentially refined meshes $\mathcal{M}_i^m$, with $i\geq 1$ 
indexing the hierarchy and $m$ indexing the subintervals or time steps, respectively. The 
initial mesh $\mathcal{M}_0^m$ is identical for each time step $m$, i.e.\ $M_0^i = 
M_0^j$ for all $i,j = \{1,\ldots, M\}$. The corresponding finite element spaces are 
denoted by $\Vs_h^{p+1,m,i}$ (cf.\ \eqref{Def:Xkhrp}) with the additional index $i$ 
denoting the mesh hierarchy. We calculate the cell- and step-wise contributions to the 
linearized error representation \eqref{Eq:DefEta} and \eqref{Thm:CellErr1}, 
respectively, by means of 
\begin{equation}
\label{Eq:DefEtaKM}
\begin{split}
&\eta_K^m = \int\limits_{I_m}\langle\mathcal{R}(u_{h}^{m,i}), z_H^{m,i}-\mathcal I_h 
z_H^{m,i}\rangle_K  - \delta_K\langle\mathcal{R}(u_{h}^{m,i}), \vec b\cdot\nabla\mathcal 
I_h z_H^{m,i}\rangle_K\\[1ex]
&\relph{\eta} \quad - \langle\mathcal{E}(u_{h}^{m,i}), z_H^{m,i}-\mathcal I_h 
z_H^{m,i}\rangle_{\partial K} \ud t  - \left\langle\left[u_{h}^{i}\right]_{m-1}, 
z_H^{m-1,i} - \mathcal I_h z_H^{m-1,i} \right\rangle_{K} \\[1ex]
&\relph{\eta} \quad + \delta_K\left\langle\left[u_{h}^i\right]_{m-1}, \vec b\cdot\nabla 
I_h z_H^{m-1,i} \right\rangle_K\,,
\end{split}\end{equation}
where the cell and edge residuals are given in \eref{eq:42} and \eref{eq:43}, 
respectively, and $\left[u_h^i\right]_0:= u_h^{1,i} - u_0$.  By $I_h 
z_H^{m,i}\in\Vs_h^{p,m,i}$ we denote the linear interpolation of the higher 
order approximation $z_H^{m,i}\in \Vs_h^{p+1,m,i}$. The integrals over the time intervals 
$I_m$ are approximated by an appropriate quadrature rule depending on the polynomial 
degree of the time discretization. For a discussion of appropriate mesh refinement 
strategies we refer to, e.g., \cite{BR03}. Details about the refinement strategy in 
time and in space that we use for the computations that are presented in Section 
\ref{Sec:NumStud} can be found in \cite{S14}. In the steady case an approach 
that is analogous to \eqref{Eq:DefEtaKM} is used. 

For the numerical computations of Section \ref{Sec:NumStud} we used the lowest 
order variants of the discretization schemes for the approximation of the primal and 
the dual solution. The discrete primal problem \eqref{Eq:dGcGs} is thus solved in 
the function space $\mathcal{X}_{kh}^{0,1}$. Up to a quadrature error in the 
right-hand side term $f$ the scheme is then algebraically equivalent to a backward Euler 
scheme in time with piecewise linear polynomials in space; cf.\ \cite{SV08} and the 
reference therein. The adjoint problem is then solved by the SUPG stabilized counterpart 
of the scheme \eqref{Eq:cGPcG} and yields a discrete solution $z_H^{i}\in 
\overline{\mathcal X}_{kh}^{\,1,2}$ being continuous and piecewise linear in time 
and continuous 
and piecewise quadratic in space. Up to a quadrature error in the right 
hand side this scheme is algebraically equivalent to the Crank-Nicolson approach 
\cite{BGR10,SV08}.    
  
For measuring the accuracy of the error estimators we will use the effectivity index 
\begin{equation}
\label{Eq:Ieff}
\mathcal I_{\mathrm eff} = \left|\frac{\eta}{\J(u)-\J(u_{kh})}\right|\,,
\end{equation}
as the ratio of the estimated error $\eta$ of \eqref{Eq:DefEta} over the exact error. 
Desirably, $\mathcal I_{\mathrm eff}$ should be close to one. In the steady case 
$\mathcal I_{\mathrm eff}$ is defined analogously with $u_{kh}$ being substituted by 
$u_{h}$.

\section{Numerical studies}
\label{Sec:NumStud}

In this section we illustrate and investigate the performance properties of the proposed 
approach of applying the Dual Weighted Residual method to stabilized finite element 
approximations of convection-dominated problems. 

\textbf{Example 1 (Hump with changing height).} As a test setting we study the moving 
hump problem that has been used in several works \cite{JS08,BS12,AM12} before as a 
benchmark problem for approximation schemes to convection-dominated equations. We 
consider problem \eqref{Eq:pro} with the prescribed solution 
\begin{equation}
\label{Eq:movghump}
\begin{split}
u(\vec x, t) &= \frac{16}{\pi}\operatorname{sin}(\pi t)x_1(1-x_1)x_2(1-x_2)\\
&\relph{=}\cdot \left\{\frac{\pi}{2} + 
\operatorname{arctan}\left(2\eps^{-\frac{1}{2}}(z_0^2 - (x_1-x_1^0)^2 - (x_2-x_2^0)^2) 
\right) \right\}\,,
\end{split}
\end{equation}
where $\Omega\times I := (0,1)^2\times(0,0.5]$ and $z_0 = 0.25$, $x_1^0 = x_2^0 = 0.5$. 
For the final time $T=0.5$ the hump reaches its maximum height. We choose the 
parameter $\eps = 10^{-6}$, $\vec b = (2,3)^\top$ and $\alpha = 1.0$. 
For the solution \eqref{Eq:movghump} the right-hand side function $f$ is 
calculated from the partial differential equation. Boundary and initial 
conditions are given by the exact solution. Our target quantity is chosen as 
\begin{equation}
\label{Eq:HumpTarget}
\J(u) = \int_\Omega u(\vec x, T) \ud \vec x\,. 
\end{equation}
We measure the spurious oscillations of the solution in the layer around the hump by  
\begin{equation}
\label{Eq:MesSpuOsc}
\mathrm{var}(t):= \underset{\vec x\in\Omega}{\mathrm{max}} \ u_{kh}(\vec x, t) - 
\underset{\vec x\in\Omega}{\mathrm{min}} \ u_{kh}(\vec x, t)\,,
\end{equation}
where the maximum and minimum are taken only in the vertices of the mesh cells. The 
exact value for the function $u$ of \eqref{Eq:movghump} at $t=0.5$ is var(0.5) = 
0.997453575; cf.\ \cite{JS08}. 

Our SUPG-stabilized discretization scheme \eqref{Eq:dGcGs} for the primal problem is 
applied with the lowest order parameter choice which amounts to $r=0$ and $p=1$. The 
dG(0) variational time discretization thus coincides with the backward Euler approach. 
For the discretization in space piecewise polynomials of first order degree are thus 
chosen. According to our derivation in Section \ref{Sec:dwr} we use a higher order 
approach with $r=1$ and $p=2$ for the discretization of the adjoint problem . 

\begin{figure}[ht!]
\centering{
\begin{tabular}{c}
\includegraphics[width=4.5cm]{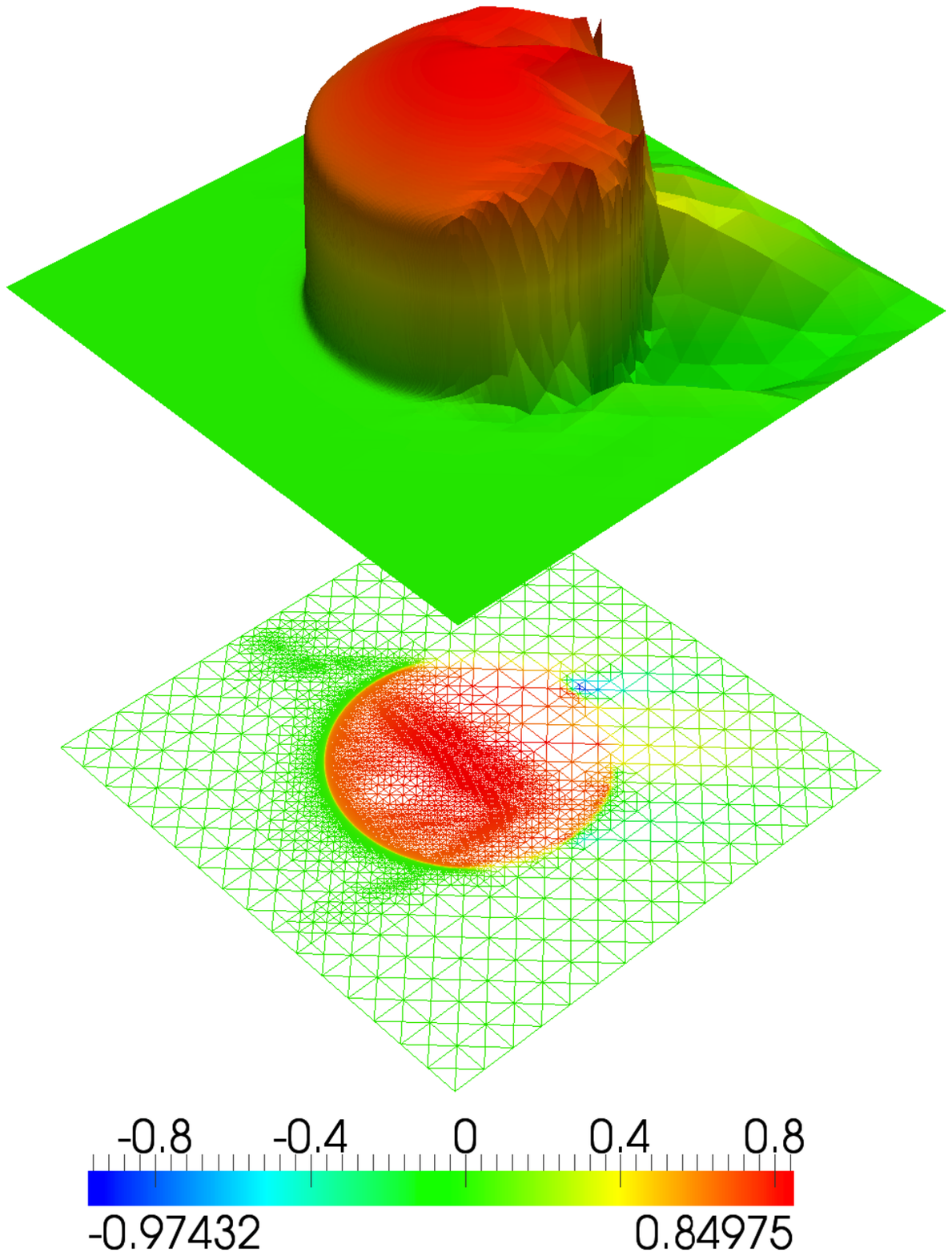}
$t = 0.25$
\end{tabular}\quad 
\begin{tabular}{c}
\includegraphics[width=4.5cm]{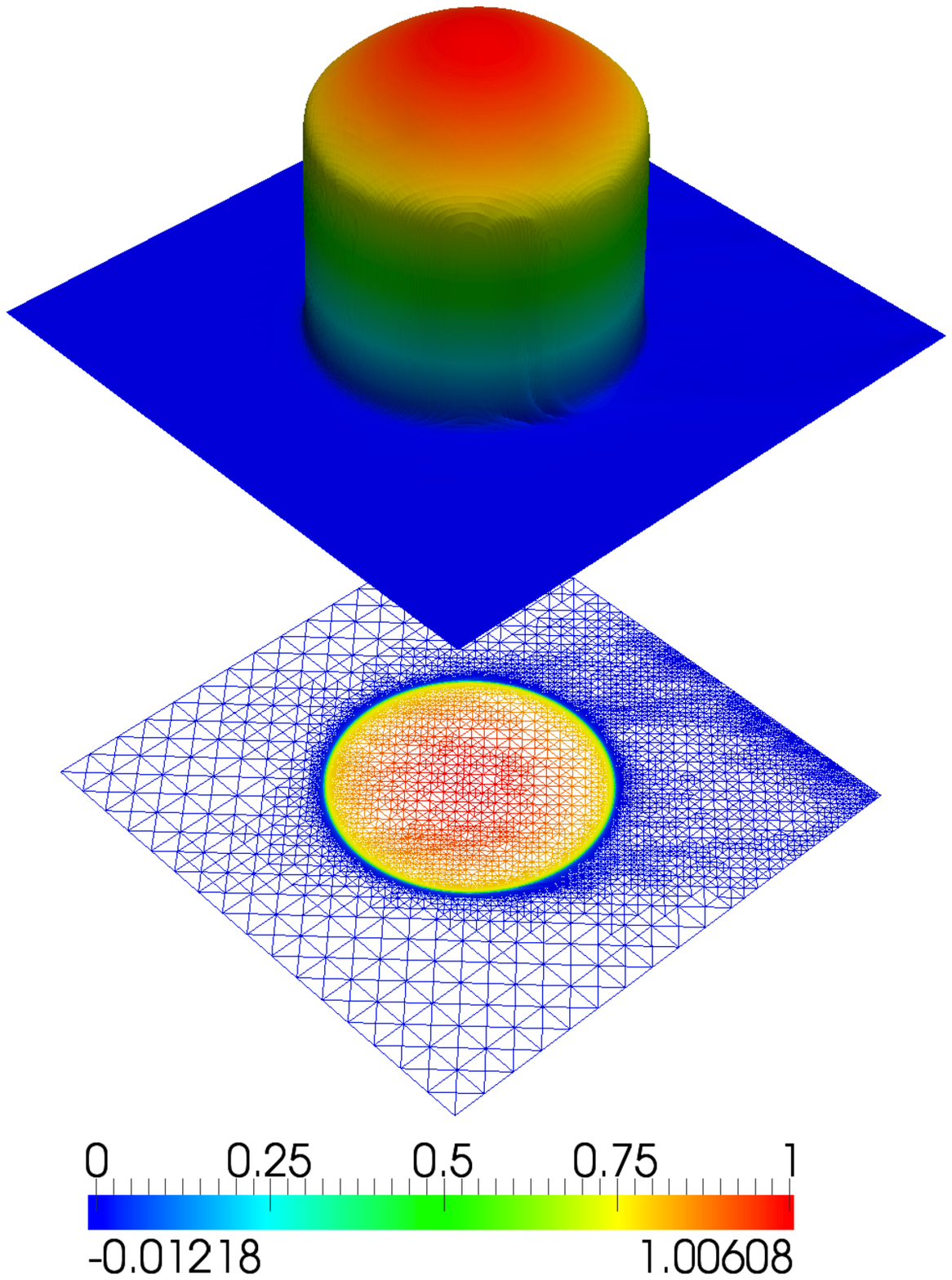}
$t = 0.5$
\end{tabular}}
\caption{Solution and adaptive mesh for goal quantity \eqref{Eq:HumpTarget} at 
time $t=0.25$ (\emph{left}) and at time $t=0.5$ (\emph{right}) for Example 1.}
\label{Fig:HumpSol}
\end{figure}

In Figure \ref{Fig:HumpSol} we visualize our computed solution profiles for the time 
points $t=0.25$ and $t=0.5$ after 16 DWR iterations on the whole time interval $(0,T]$. 
For $t=0.25$ the solution is still strongly perturbed in the backward part of the hump's 
layer and behind the hump in the direction of the flow field $\vec b$. The mesh is coarse 
in that part of the domain. Such a behaviour is admissible since our target functional 
aims to control the solution profile at the final time point $T=0.5$ only. For $T=0.5$ an 
almost perfect solution profile is obtained and the finite element mesh cells are 
concentrated on the backward face of the hump. We note that the spurious oscillations 
behind the hump, that were obtained by different classes of approximation schemes in 
\cite{JS08}, do not arise here. They are strongly reduced and almost completely eliminated 
by the adaptive algorithm. In Figure \ref{Fig:HumpCompData_1} the magnitude of the 
adaptively chosen time steps is presented. The first time steps are chosen relatively 
large whereas the time step sizes close to the final time point $T=0.5$ become much 
smaller. Even though large time step sizes and also large spatial mesh sizes in the 
crucial regions are used in the first time steps, leading to crude approximations in the 
initial phase as shown in the left plot of Figure \ref{Fig:HumpSol}, the algorithm is 
capable to provide the desired approximation quality in the target quantity 
\eqref{Eq:HumpTarget} that is local in time and controls the solution profile at the time 
final time point only. A high approximation quality in the target quantity is thus 
obtained with very economical meshes. 

In Table \ref{Tab:HumpCompData_3} we monitor the convergence behaviour of the 
DWR iterations in terms of the effectivity index \eqref{Eq:Ieff} for the target 
functional \eqref{Eq:HumpTarget}. The degrees of freedom at the final time point $T=0.5$ 
are given together with the corresponding values of the effectivity index and 
$\mathrm{var}(0.5)$. For an increasing number of DWR iterations with space-time 
mesh adaptions the effectivity index is very close to one indicating an excellent 
approximation of the goal quantity by the DWR approach applied to the stabilized 
approximation of \eqref{Eq:pro}. Further, the given numbers for $\mathrm{var}(0.5)$ show 
that the spurious oscillations in the layer around the hump are also reduced by the DWR 
iterations and the space-time grid adaption process. This might be a consequence 
of the global character in space of the target functional \eqref{Eq:HumpTarget}. Finally 
in Table \ref{Tab:HumpCompData_2} we compare the values for $\mathrm{var}(0.5)$ that we 
computed by our adaptive approach with some reference values that were obtained by other 
research groups and published in the literature. The calculations of all other groups 
were done on uniform meshes. For comparison purposes our adaptive simulations were 
run in such a way that either the number of degrees of freedom 
or the calculated value $\mathrm{var}(0.5)$ coincides approximately with the 
given reference values of the literature. The presented numbers impressively 
illustrate the superiority of the adaptive computations. 

\begin{minipage}{\linewidth}
\begin{minipage}{5cm}
\begin{figure}[H]
\includegraphics[width=5cm]{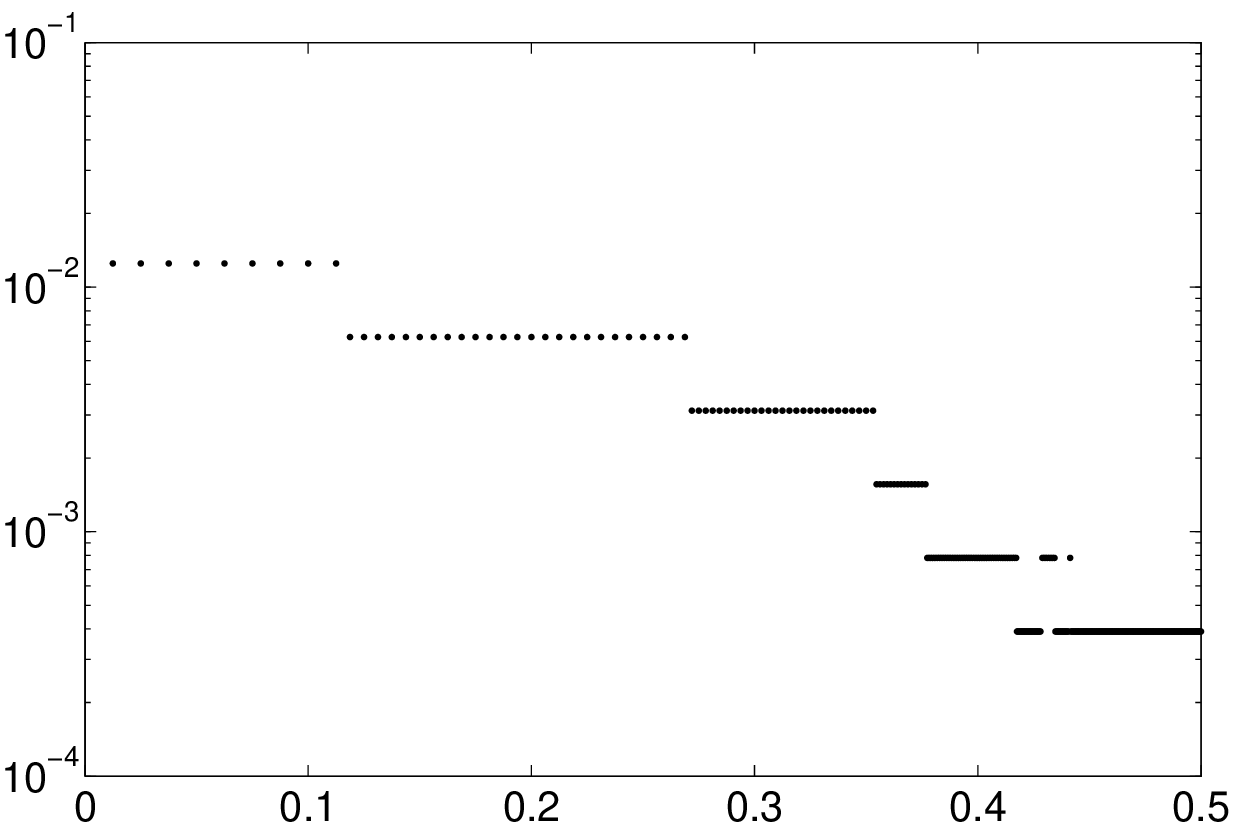}
\caption{Time step sizes over $(0,T]$ for Example 1.}
\label{Fig:HumpCompData_1}
\end{figure}
\end{minipage}\hfill 
\begin{minipage}{7.4cm}
\begin{table}[H]
\footnotesize \centering
$\begin{array}{ccc|ccc}
\toprule
\toprule
\mbox{dofs} & \mathrm{var}(0.5) & \mathcal I_{\mathrm eff} & \mbox{dofs} 
& \mathrm{var}(0.5) & \mathcal I_{\mathrm eff}\\
\midrule
  545 & 9.392 & 1.03 &  6425 & 1.200 & 1.02 \\
  833 & 9.432 & 1.17 &  7283 & 1.115 & 1.00 \\
 1267 & 5.852 & 1.75 &  8939 & 1.076 & 0.96 \\
 1535 & 6.761 & 9.13 & 11954 & 1.079 & 0.97 \\
 2325 & 4.328 & 2.82 & 20757 & 1.050 & 0.99 \\
 2522 & 3.270 & 1.08 & 21455 & 1.030 & 1.03 \\
 3031 & 2.105 & 1.47 & 25492 & 1.022 & 1.05 \\
 3937 & 1.568 & 0.94 & 36167 & 1.018 & 1.07 \\
\bottomrule
\end{array}$
\caption{Convergence statistics for Example 1.}
\label{Tab:HumpCompData_3}
\end{table}
\end{minipage}
\end{minipage}

\begin{table}[ht!]
\centering{
\begin{minipage}[t]{0.7\textwidth}
\footnotesize
\centering
\begin{tabular}{>{$}l<{$}>{$}c<{$}>{$}c<{$}>{$}r<{$}>{$}r<{$}}
\toprule
\toprule
\mbox{Method} & \mbox{Reference} & \mathrm{var}(0.5) & \multicolumn{1}{c}{\mbox{dofs}} 
& \multicolumn{1}{c}{$k$}\\
\midrule
\mbox{SUPG}     & \cite{JS08} & 1.3835 & 16641 & 10^{-3}\\
\mbox{LPS}      & \cite{JS08} & 1.2007 & 32768 & 10^{-3}\\
\mbox{SUPG}     & \cite{BS12} & 1.2504 & 33025 & 2\cdot 10^{-3}\\
\mbox{SUPG/SC}  & \cite{BS12} & 1.1946 & 33025 & 2\cdot 10^{-3}\\
\mathrm{LPS/cGP}(1) & \cite{AM12} & 1.0408 & 33025 & 10^{-3}\\
\mathrm{LPS/dG}(1)  & \cite{AM12} & 1.0408 & 33025 & 10^{-3}\\
\mbox{SUPG/DWR} & \mbox{ this work} & 1.0790 & 10900 & 3.1\cdot 10^{-3}\\
\mbox{SUPG/DWR} & \mbox{ this work} & 1.0179 & 35931& 1.5\cdot 10^{-3}\\
\bottomrule
\end{tabular}
\caption{Reference values of the literature for Example 1.}
\label{Tab:HumpCompData_2}
\end{minipage}
}
\end{table} 

\textbf{Example 2 (Point-value error control).} In this example we illustrate the 
application of our approach to a target functional that provides a spatially local error 
control in a sharp layer. Thereby we evaluate the potential of our approach to capture 
sharp layers and fronts with high accuracy. This is a challenging task and of utmost 
interest for convection-dominated problems. Since the interaction of the goal-oriented 
error control mechanism with the discretization in space is especially focused here, we 
restrict ourselves to the stationary case for simplicity. As a benchmark problem we use 
an 
adaptation of \cite[Example 4.2]{LR06}. We consider problem \eqref{Eq:pros} with $\Omega 
= (0,1)^2$, $\alpha = 1.0$, $\vec b = \frac{1}{\sqrt{5}}(1,2)^{\top}$, $\varepsilon = 
10^{-6}$ and nonlinear reaction term $r(u) = u^2$. We choose the right-hand side $f$ such 
that 
\begin{equation}
\label{Eq:SolExp2}
u(\vec x) = \frac{1}{2}\left(1-\operatorname{tanh}\frac{2x_1 - 
x_2-0.25}{\sqrt{5\eps}}\right)
\end{equation}
is the analytical solution of \eqref{Eq:pros}. The Dirichlet boundary condition is given 
by the exact solution. The solution is characterized by an interior layer of thickness 
$\mathcal O(\sqrt{\varepsilon}\vert \operatorname{ln} \varepsilon\vert)$. We study the 
following target functionals 
\begin{align*}
&\J_{\Ls^2}(u) = \frac{1}{\Vert e\Vert_{\Ls^2(\Omega)}}\langle 
e,u\rangle_{\Omega}\,, \qquad \J_1(u) = \int_\Omega u\ud\vec x \qquad \mbox{ and } 
&\J_2(u) = u(\vec x_e)\,,
\end{align*}
where $e:=u-u_h$ and with a user-prescribed control point $\vec x_e 
= \left(\frac{3}{16},\frac{1}{8}\right)$ that is located in the interior of the layer. In 
our computations we regularize the functional $\J_2(u)$ by 
\begin{align*}
\J_{r}(u) = \frac{1}{\vert B_r\vert}\int_{B_r} u(\vec x)\ud\vec x\,,
\end{align*}
where the ball $B_r$ is defined by $B_r = \left\{\vec x\in\Omega \,\vert\, \vert\vec 
x-\vec x_e\vert < r \right\}$ with a small radius $r$.

\begin{figure}[ht!]
\centering
\begin{tabular}{cc}
\includegraphics[width=4cm]{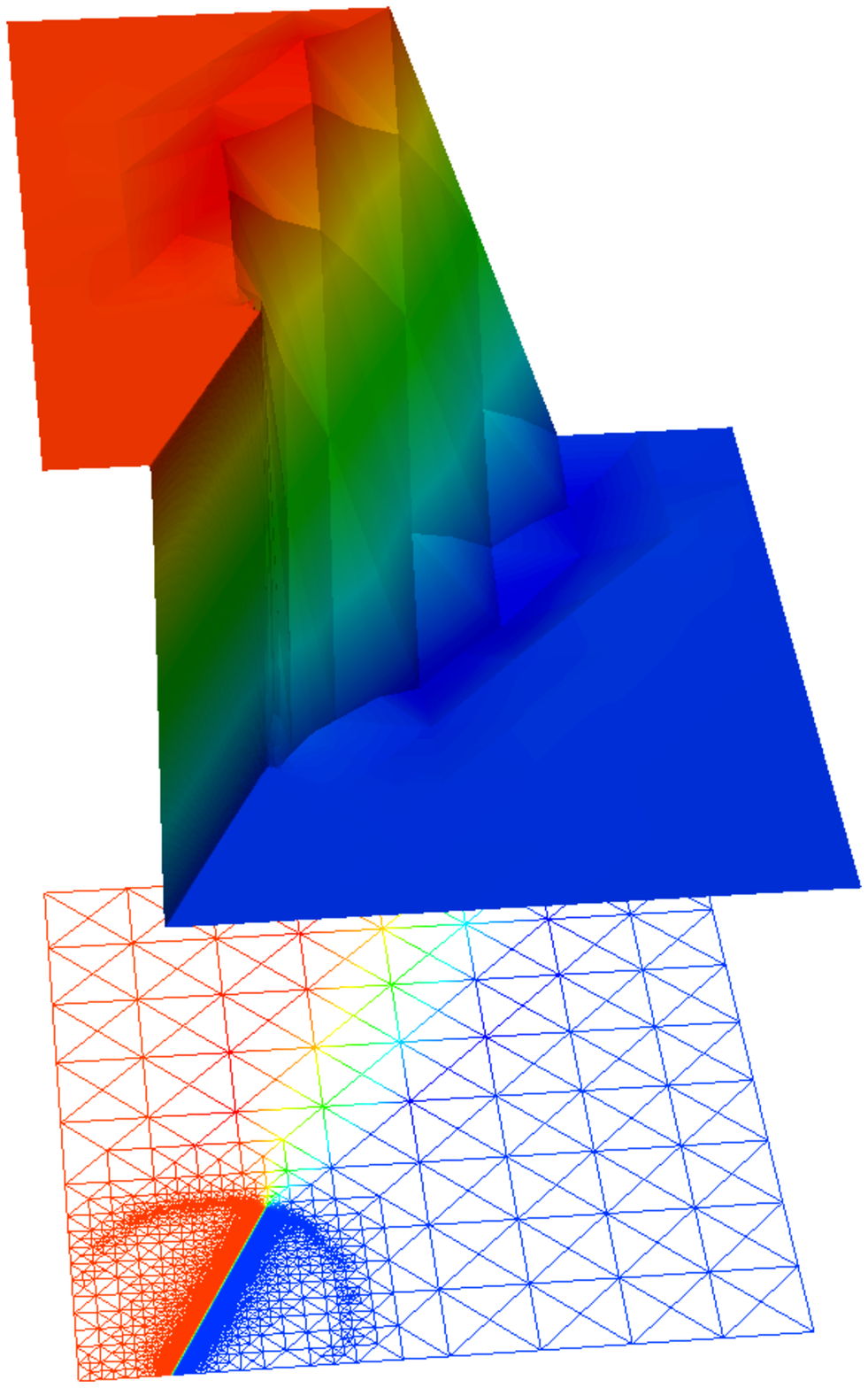}&
\includegraphics[width=4cm]{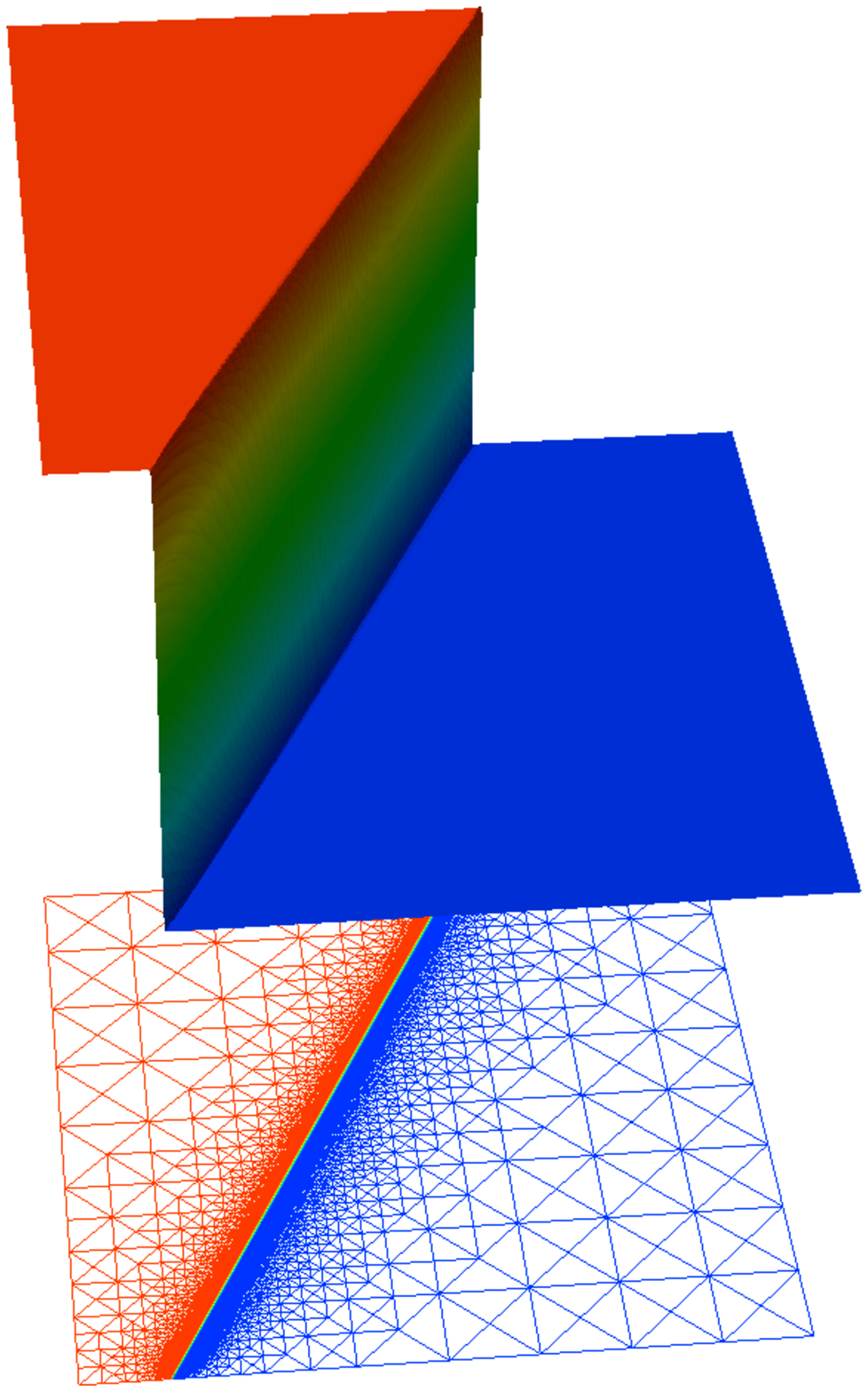}\\
Point-value error control. & Global error control.
\end{tabular}
\caption{Point-value error control (\emph{left}) and global error control 
(\emph{right}) by the DWR approach for Example 2.}
\label{Fig:TanhSol}
\end{figure}

In Figure \ref{Fig:TanhConv} and Table \ref{Tab:TanhIeff} we summarize 
the convergence behavior of proposed DWR approach to the stabilized approximation scheme 
\eqref{eq:08}. We note that $\J_1(\cdot)$ provides the traditional global 
$\mathcal L^2$-error 
control and is considered for reference purposes. For the target functionals 
$\J_{\Ls^2}(\cdot)$ and $\J_1(\cdot)$ the effectivity indices converge to one for an 
increasing number of degrees of freedom. For the challenging point-value error control 
of $\J_r(\cdot)$ the effectivity index is also very close to one which is in good 
agreement with effectivity indices for point-value error control that are given in other 
works of the literature; cf.\ \cite[p.\ 45]{BR03} for the pure Poisson problem. In Figure 
\ref{Fig:TanhSol} we visualize the computed solution profiles and adaptive meshes for an 
error control based on the local target functional $\J_r$ and the global target 
functional 
$\J_1$, respectively. This example nicely brings out the potential of the DWR approach. 
For the point-value error control the mesh cells are located around the specified point 
of 
interest. Even though a crude approximation of the sharp interface is obtained away from 
the specified control point, in its neighborhood an excellent approximation of the sharp 
layer is ensured by the approach. A very economical mesh along with a high quality in the 
computation of the user-specified goal quantity is thus obtained. The global error 
control 
of $J_1$ provides a good approximation of the solution in the whole domain by adjusting 
the mesh along the layer.    

\begin{minipage}{\linewidth}
\begin{minipage}{5.3cm}
\begin{figure}[H]
\centering
\includegraphics[width=5cm]{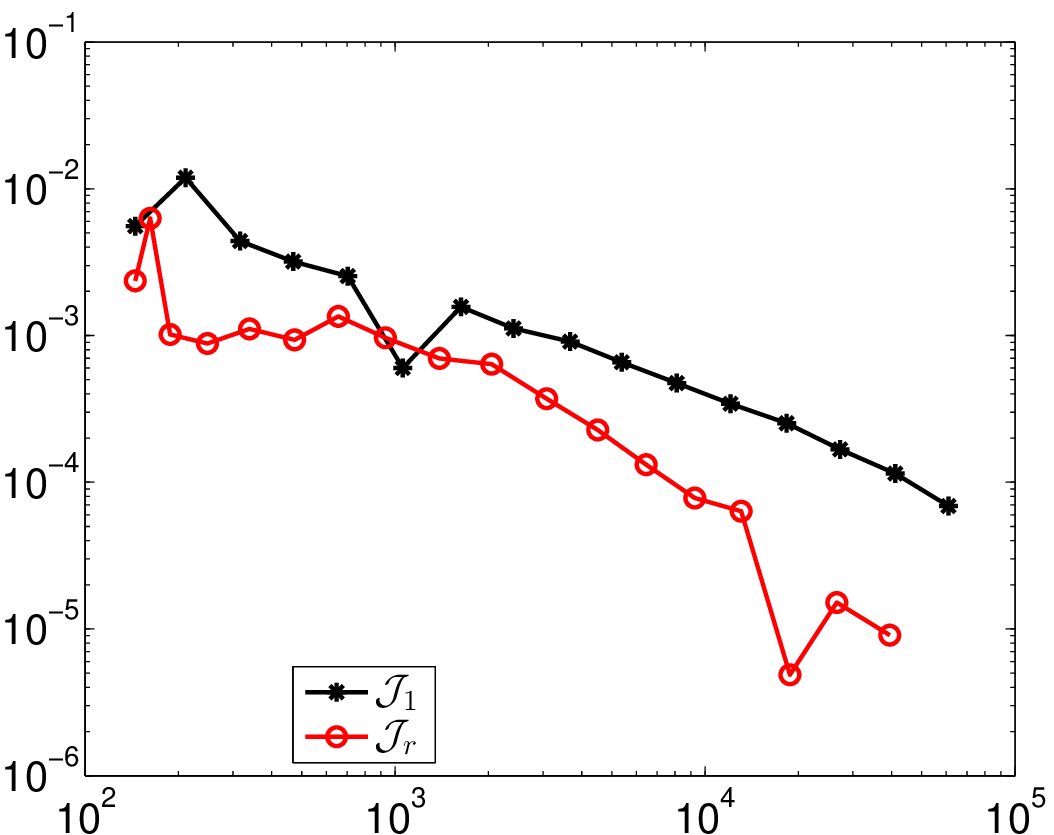}
\caption{$\J_1(u-u_h)$ and $\J_r(u-u_h)$ over degrees of freedom for Example 2.}
\label{Fig:TanhConv}
\end{figure}
\end{minipage}\hfill 
\begin{minipage}{7.4cm}
\begin{table}[H]
\footnotesize \centering
$\begin{array}{cc|cc|cc}
\toprule
\toprule
\multicolumn{2}{c|}{\J_{\Ls^2}} & \multicolumn{2}{c|}{\J_1} & \multicolumn{2}{c}{\J_r}\\
\midrule
\mbox{dofs} & \mathcal I_{\mathrm eff} & \mbox{dofs} & \mathcal I_{\mathrm eff} 
& \mathcal I_{\mathrm eff} & \mbox{dofs} \\
\midrule
\midrule
 3756 & 0.74 &  5383 & 0.45 & 0.03  & 4505 \\
 5903 & 0.83 &  8105 & 0.44 & 0.07  & 6458 \\
 9059 & 0.87 & 12081 & 0.45 & 0.14 &  9268 \\
14373 & 0.94 & 18321 & 0.57 & 0.22 & 13079 \\
22834 & 0.95 & 27276 & 0.71 & 3.42 & 18794 \\
37555 & 0.97 & 41073 & 0.76 & 1.25 & 26619 \\
62119 & 0.98 & 60957 & 0.83 & 2.27 & 39447 \\
\bottomrule
\end{array}$
\caption{Effectivity indices for the target quantities $\J_{\Ls^2}$\,, $\J_1$ 
and $\J_r$ for Example 2.}
\label{Tab:TanhIeff}
\end{table}
\end{minipage}
\end{minipage}

\begin{figure}[ht!]
\centering
\includegraphics[width=6.50cm]{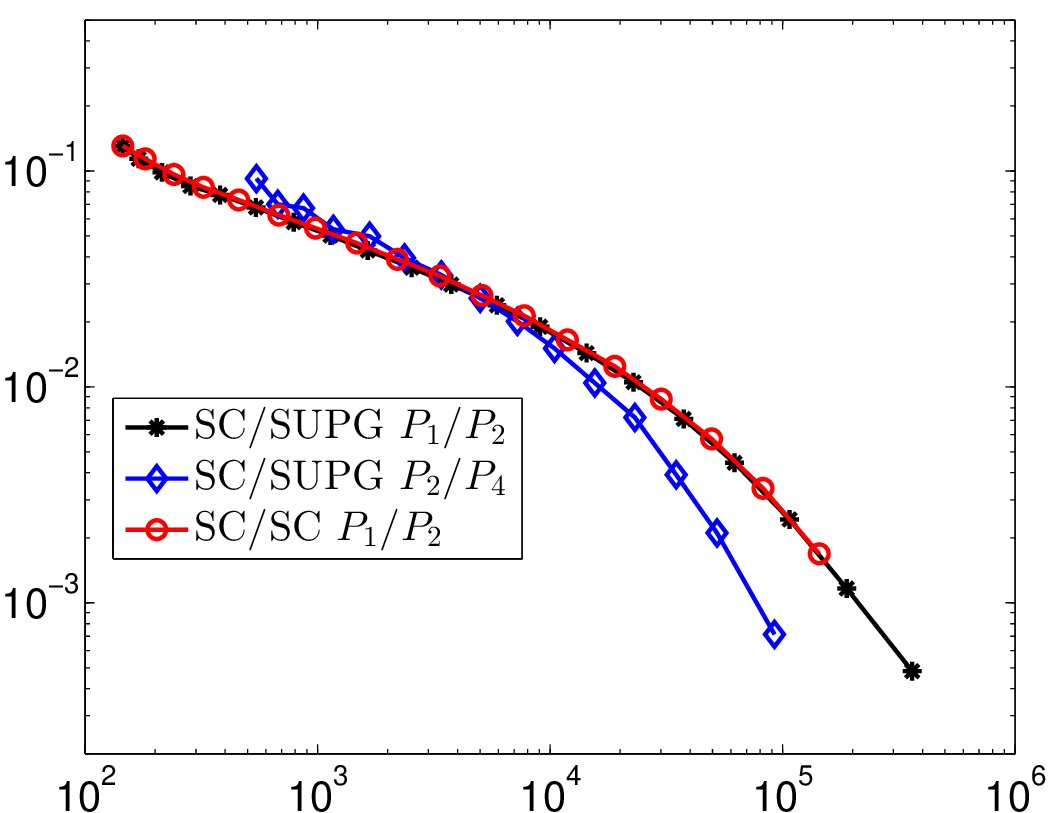}
\includegraphics[width=6.4cm]{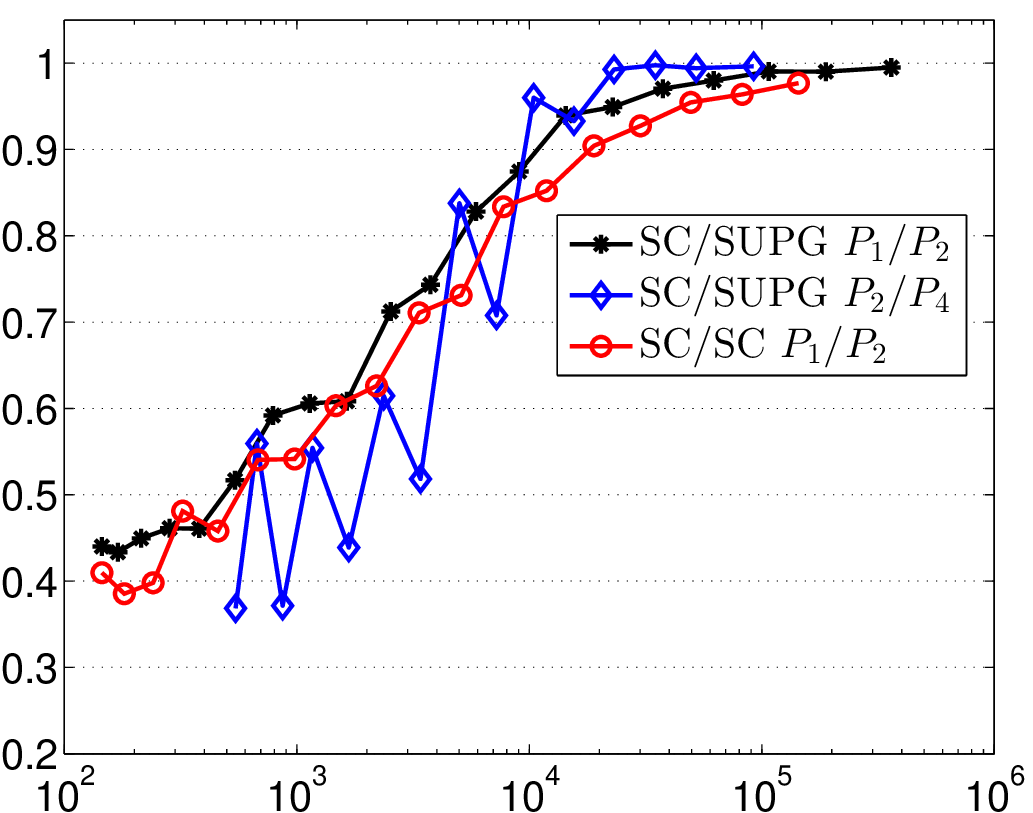}
\caption{Errors $\J_{\Ls^2}(u-u_h)$ (\emph{left}) and effectivity indices 
$\mathcal{I}_{\mathrm{eff}}$ (\emph{right}) over degrees of freedom for 
Example 3 with stabilization techniques for the primal\slash dual problem and 
polynomial degrees for the approximation of the primal\slash dual problem -- SC 
indicates the application of additional shock--capturing stabilization, SUPG 
means SUPG stabilization only).}
\label{Fig:ImpDualApprox}
\end{figure}

\textbf{Example 3 (Impact of approximation of the weights).} In the last example we 
briefly study the impact of the approximation of the dual solution $z$ on the quality of 
the overall error control process. For brevity we consider again the test problem of 
Example 2 with the prescribed solution \eqref{Eq:SolExp2}. In our first study we use a 
piecewise linear approximation in $\mathcal V_h^{1}$ with SUPG and additional 
shock-capturing stabilization for the primal problem. For the corresponding adjoint 
problem a piecewise quadratic finite approximation in $\mathcal V_h^{2}$ with SUPG 
stabilization and with and without additional shock-capturing stabilization is used. The 
left plot of Figure \ref{Fig:ImpDualApprox} shows that the additional shock-capturing 
stabilization of the dual problem yields no further improvement in the accuracy of the 
approximation. This is advantageous since the adjoint problem by itself is always 
a linear one. Applying  shock-capturing stabilization introduces an artificial 
nonlinearity and requires (nonlinear) iterations for solving the arising algebraic 
system. The left plot of Figure \ref{Fig:ImpDualApprox} that argues that using only SUPG 
stabilization for the dual problem and thereby keeping its linear character is 
sufficient for the proposed DWR apporach. We note that the positive impact of additional 
shock-capturing stabilization in the numerical approximation of convection-dominated 
problems has been well understood and analyzed numerically; cf.\ \cite{BS12,JS08} and the 
references therein. 

Further, the left plot of Figure \ref{Fig:ImpDualApprox} shows the gain in accuracy if a 
higher order approach is used. Here we combined a stabilized piecewise quadratic 
approximation in $\mathcal V_h^{2}$ of the primal problem with a stabilized approximation 
in $\mathcal V_h^{4}$ with piecewise polynomials of fourth order of the adjoint problem. 
An approximation of the adjoint problem with piecewise polynomials of third order 
did not provide sufficient accuracy and did not yield a convergence behaviour or an error 
reduction, respectively, similarly to the one that is shown in the left plot of Figure 
\ref{Fig:ImpDualApprox}. This observation underlines the necessity of the proper 
approximation of the adjoint problem within the DWR framework. For non 
convection-dominated problems the process might be not that much sensitive as in our 
studies for problems with strong layers and sharp fronts. In the right plot of Figure 
\ref{Fig:ImpDualApprox} the corresponding values of the effectivity index are visualized.

\section{Summary}

In this work we developed an adaptive approach for stabilized finite element 
approximations of convection-dominated problems. It is based on the dual 
weighted residual method for goal-oriented a posteriori error control. 
A \emph{first dualize and then stabilize} philosophy was applied for combining 
the Dual Weighted Residual method with the stabilization of the finite element 
techniques. In contrast to other works of the literature we used a higher order 
approximation of the adjoint problem instead of a higher order interpolation of 
a lower order approximation of the dual solution. Thereby we aim to eliminate 
sources of inaccuracies in regions with layers and close to sharp fronts. 
In numerical experiments we could prove that spurious oscillations that 
typically arise in numerical approximations of convection-dominated problems 
could be reduced significantly. Effectivity indices very close to one were 
obtained for the user-specified target quantities. The presented approach offers 
large potential for combining goal-oriented error control and selfadaptivity 
with stabilized finite element methods in the approximation of 
convection-dominated transport. The application of the approach to more 
sophisticated problems, like Navier--Stokes problems, is our work for the future. 
Moreover, the efficient computation of the higher order approximation to the 
adjoint problem offers potential for optimization. 
This will also be our work for the future.

\appendix
\numberwithin{equation}{section}
\section{Appendix}

For the sake of completeness we provide the proof of Theorem \ref{Thm:L}.

\begin{proof} 
We let $e = x_c -x_d$. By the fundamental theorem of calculus it holds that 
\[
\mathcal L(x_c) - \mathcal L(x_d) = \int_0^1 \mathcal L'(x_d+s e)(e)\ud s\,.
\]
Approximating the integral by the trapezoidal rule yields that 
\begin{equation}
\label{Eq:EL10}
\mathcal L(x_c) - \mathcal L(x_d) = \frac{1}{2} \mathcal L'(x_d)(x_c -x_d) + \frac{1}{2} 
\mathcal L'(x_c)(x_c -x_d) + \mathcal R
\end{equation}
with $\mathcal R$ being defined by \eqref{Def:L4}. By the supposed stationarity 
of $\mathcal{L}$ in $x_c$ along with the assumption \eqref{Def:L1} the second 
of terms on the right-hand side of \eqref{Eq:EL10} vanishes. Together with 
eq.\ \eqref{Def:L0} we then 
get that 
\begin{align*}
\mathcal L(x_c) - \mathcal L(x_d) & = \frac{1}{2} \mathcal L'(x_d) (x_c- y_d) +  
\frac{1}{2} \mathcal L'(x_d) (y_d - x_d) + \mathcal R\\[1ex]
& = \frac{1}{2} \mathcal L'(x_d) (x_c- y_d) + \frac{1}{2} \mathcal 
D(x_d)(y_d-x_d) + 
\mathcal R
\end{align*}
for all $y_d\in X_d$. This completes the proof of Theorem \ref{Thm:L}.
\end{proof}

\end{document}